\documentclass[]{article}

\usepackage{lineno,hyperref}
\modulolinenumbers[5]

\usepackage[textwidth=1.5in,textsize=small]{todonotes}
\usepackage[margin=2.54cm]{geometry}
\usepackage{graphicx}
\usepackage{subfig}
\usepackage{mathrsfs}
\usepackage{tikz}
\usetikzlibrary{matrix,backgrounds,shapes}
\usetikzlibrary{fit}
\usepackage{xspace}

\usepackage{amssymb}
\usepackage{amsmath}
\usepackage{amsthm}
\usepackage{mathtools}
\usepackage{bm}
\definecolor{mayablue}{rgb}{0.45, 0.76, 0.98}
\usepackage{epstopdf}
\usepackage{stmaryrd}
\usepackage{lmodern}
\usepackage{enumitem}
\usepackage{amsfonts}
\usepackage{booktabs}
\usepackage{algorithm2e}

\newcommand{\textred}[1]{\textcolor{black}{#1}}
\newcommand{\textblue}[1]{\textcolor{black}{#1}}

\newcommand{\fnc}[1]{\ensuremath{\mathcal{#1}}}
\newcommand{\mat}[1]{\ensuremath{\mathsf{#1}}}
\newcommand{\Tr}[0]{^{\textrm{T}}}
\newcommand{\mr}[1]{\ensuremath{\mathrm{#1}}}

\newcommand{\Nxi}[0]{N_{\xi}}
\newcommand{\Neta}[0]{N_{\eta}}

\newcommand{\nmin}[1]{N^{*}\left(#1\right)}
\newcommand{\Pk}[0]{\fnc{P}_{k}}
\newcommand{\bpk}[0]{\bm{p}_{k}}

\newcommand{\bpkL}[0]{\bm{p}_{k,L}}
\newcommand{\bpkGamL}[0]{\bm{p}_{k,\hat{\Gamma},L}}

\newcommand{\bpktR}[0]{\tilde{\bm{p}}_{k,R}}

\newcommand{\dxibpk}[0]{\frac{\partial \bpk}{\partial \xi} }
\newcommand{\etakL}[0]{\bm{\eta}_{L}^{k}}

\newcommand{\etakR}[0]{\bm{\eta}_{R}^{k}}

\newcommand{\etakGam}[0]{\bm{\eta}_{\Gamma}^{k}}
\newcommand{\M}[0]{\mat{H}}

\newcommand{\Dxi}[0]{\mat{D}_{\xi}}
\newcommand{\Deta}[0]{\mat{D}_{\eta}}
\newcommand{\Qxi}[0]{\mat{Q}_{\xi}}
\newcommand{\Qeta}[0]{\mat{Q}_{\eta}}
\newcommand{\Sxi}[0]{\mat{S}_{\xi}}
\newcommand{\Seta}[0]{\mat{S}_{\eta}}
\newcommand{\Exi}[0]{\mat{E}_{\xi}}
\newcommand{\Eeta}[0]{\mat{E}_{\eta}}
\newcommand{\nxi}[0]{n_{\xi}}
\newcommand{\nxiL}[0]{n_{\xi,L}}
\newcommand{\neta}[0]{n_{\eta}}

\newcommand{\ML}[0]{\mat{H}_{L}}
\newcommand{\DxiL}[0]{\mat{D}_{\xi_{L}}}

\newcommand{\SxiL}[0]{\mat{S}_{\xi_{L}}}

\newcommand{\MR}[0]{\mat{H}_{R}}
\newcommand{\DxiR}[0]{\mat{D}_{\xi_{R}}}

\newcommand{\etaL}[0]{\bm\eta_{L}}
\newcommand{\etaR}[0]{\bm\eta_{R}}
\newcommand{\etaGam}[0]{\bm\eta_{\Gamma}}
\newcommand{\etaLk}[0]{\bm\eta^k_{L}}
\newcommand{\etaRk}[0]{\bm\eta^k_{R}}
\newcommand{\etaGamk}[0]{\bm\eta^k_{\Gamma}}

\newcommand{\Dxione}[0]{\mat{D}_{\xi}^{(1\textrm{D})}}
\newcommand{\DxioneL}[0]{\mat{D}_{\xi_L}^{(1\textrm{D})}}
\newcommand{\Detaone}[0]{\mat{D}_{\eta}^{(1\textrm{D})}}

\newcommand{\Qxione}[0]{\mat{Q}_{\xi}^{(1\textrm{D})}}
\newcommand{\Qetaone}[0]{\mat{Q}_{\eta}^{(1\textrm{D})}}
\newcommand{\Sxione}[0]{\mat{S}_{\xi}^{(1\textrm{D})}}
\newcommand{\Setaone}[0]{\mat{S}_{\eta}^{(1\textrm{D})}}

\newcommand{\Mxione}[0]{\mat{H}_{\xi}^{(1\textrm{D})}}
\newcommand{\Metaone}[0]{\mat{H}_{\eta}^{(1\textrm{D})}}

\newcommand{\MxiLone}[0]{\mat{H}_{\xi_{L}}^{(1\textrm{D})}}

\newcommand{\MetaLone}[0]{\mat{H}_{\eta_{L}}^{(1\textrm{D})}}
\newcommand{\MetaRone}[0]{\mat{H}_{\eta_{R}}^{(1\textrm{D})}}

\newcommand{\exioneR}[0]{\bm{e}_{\xi_R,1}}

\newcommand{\exiNL}[0]{\bm{e}_{\xi_L,\Nxi}}

\newcommand{\Ixi}[0]{\mat{I}_{\xi}}

\newcommand{\OneL}[0]{\bm{1}_{L}}
\newcommand{\OneR}[0]{\bm{1}_{R}}
\newcommand{\Ieta}[0]{\mat{I}_{\eta}}

\newcommand{\ExioneD}[0]{\mat{E}_{\xi}^{(1\textrm{D})}}

\newcommand{\Exione}[0]{\mat{E}_{\xi,1}}
\newcommand{\ExiN}[0]{\mat{E}_{\xi,\Nxi}}
\newcommand{\Eetaone}[0]{\mat{E}_{\eta,1}}
\newcommand{\EetaN}[0]{\mat{E}_{\eta,\Neta}}
\newcommand{\Rxione}[0]{\mat{R}_{\xi,1}}
\newcommand{\RxiN}[0]{\mat{R}_{\xi,\Nxi}}
\newcommand{\Retaone}[0]{\mat{R}_{\eta,1}}
\newcommand{\RetaN}[0]{\mat{R}_{\eta,\Neta}}

\newcommand{\RxiLN}[0]{\mat{R}_{\xi_L,\Nxi}}

\newcommand{\RxiRone}[0]{\mat{R}_{\xi_R,1}}

\newcommand{\Mclas}{\mat{H}_{c}^{(1\textrm{D})}}
\newcommand{\Dclas}{\mat{D}_{c}^{(1\textrm{D})}}
\newcommand{\Qclas}{\mat{Q}_{c}^{(1\textrm{D})}}

\newcommand{\betax}[0]{\alpha}
\newcommand{\betay}[0]{\beta}
\newcommand{\nx}[0]{n_{x}}
\newcommand{\ny}[0]{n_{y}}
\newcommand{\Gamman}[0]{\Gamma^{-}}
\newcommand{\Gammap}[0]{\Gamma^{+}}

\newcommand{\lamxi}[0]{\tilde\alpha}
\newcommand{\lameta}[0]{\tilde\beta}
\newcommand{\lamxiL}[0]{\tilde\alpha_{L}}
\newcommand{\lametaL}[0]{\tilde\beta_{L}}
\newcommand{\lamxiR}[0]{\tilde\alpha_{R}}
\newcommand{\lametaR}[0]{\tilde\beta_{R}}
\newcommand{\uL}[0]{\bm{u}_{L}}
\newcommand{\uR}[0]{\bm{u}_{R}}
\newcommand{\J}[0]{\frac{1}{\fnc{J}}}
\newcommand{\JL}[0]{\frac{1}{\fnc{J}_{L}}}

\newcommand{\ExiNL}[0]{\mat{E}_{\xi_{L},\Nxi}}

\newcommand{\ExioneR}[0]{\mat{E}_{\xi_{R},1}}

\newcommand{\RxioneR}[0]{\mat{R}_{\xi_{R},1}}

\newcommand{\xik}[1]{\bm{\xi}^{#1}}
\newcommand{\xil}[0]{\xi_{L}}

\newcommand{\txil}[0]{\bm{e}_{\xi_{1}}}
\newcommand{\txir}[0]{\bm{e}_{\xi_{N}}}

\newcommand{\txione}[0]{\bm{e}_{\xi,1}}
\newcommand{\txiN}[0]{\bm{e}_{\xi,\Nxi}}
\newcommand{\tetaone}[0]{\bm{e}_{\eta,1}}
\newcommand{\tetaN}[0]{\bm{e}_{\eta,\Neta}}

\newcommand{\RL}[0]{\mat{P}_{\xi,L}}
\newcommand{\RR}[0]{\mat{P}_{\xi,R}}
\newcommand{\RLone}[0]{\mat{P}_{\xi,L}^{(1\textrm{D})}}
\newcommand{\RRone}[0]{\mat{P}_{\xi,R}^{(1\textrm{D})}}
\newcommand{\RLoneT}[0]{\mat{P}_{\xi,L}^{(1\textrm{D}),T}}
\newcommand{\RRoneT}[0]{\mat{P}_{\xi,R}^{(1\textrm{D}),T}}
\newcommand{\MGam}[0]{\mat{M}_{\hat{\Gamma}}}
\newcommand{\diag}[0]{\mr{diag}}

\newcommand{\uLN}{\bm{u}_{L,N}}
\newcommand{\uRone}{\bm{u}_{R,1}}
\newcommand{\uGamL}{\bm{u}_{L,\hat{\Gamma}}}
\newcommand{\uGamR}{\bm{u}_{R,\hat{\Gamma}}}
\newcommand{\Ohat}[0]{\hat{\Omega}}
\newcommand{\Ghat}[0]{\hat{\Gamma}}

\newcommand{\HA}[0]{\mat{H}_{L}}
\newcommand{\HB}[0]{\mat{H}_{R}}

\newcommand{\DxiA}[0]{\mat{D}_{\xi_{L}}}
\newcommand{\DxiB}[0]{\mat{D}_{\xi_{R}}}
\newcommand{\DetaA}[0]{\mat{D}_{\eta_{L}}}
\newcommand{\DetaB}[0]{\mat{D}_{\eta_{R}}}

\newcommand{\uA}[0]{\bm{u}_{L}}
\newcommand{\uB}[0]{\bm{u}_{R}}
\newcommand{\JA}[0]{\frac{1}{\fnc{J}_{L}}}
\newcommand{\JB}[0]{\frac{1}{\fnc{J}_{R}}}

\newcommand{\lamxiA}[0]{\tilde \alpha_{L}}
\newcommand{\lamxiB}[0]{\tilde\alpha_{R}}
\newcommand{\lametaA}[0]{\tilde\beta_{L}}
\newcommand{\lametaB}[0]{\tilde\beta_{R}}

\newcommand{\ignore}[1]{}

\newcommand{\fsL}[0]{\bm{f}_{L}^{*}}
\newcommand{\fsR}[0]{\bm{f}_{R}^{*}}
\newcommand{\fs}[0]{\bm{f}^{*}}

\newcommand{\fsGas}[0]{\bm{f}_{Gas}^{*}}
\newcommand{\Ltwo}[0]{L_{2}}
\theoremstyle{plain}
\newtheorem{thm}{Theorem}
\newtheoremstyle{proofpartstyle}
  {1pc} 
  {1pc} 
  {} 
  {2.5mm} 
  {\itshape} 
  {} 
  {-0.1em} 
  {} 
\theoremstyle{proofpartstyle} \newtheorem*{proofpart}{Proof of Part}
\makeatletter
\@addtoreset{proofpart}{theorem}
\makeatother

\theoremstyle{plain}
\newtheorem{definition}{Definition}
\newtheorem{rmk}{Remark}
\newtheorem{pf}{Proof}
\numberwithin{equation}{section}

\begin{document}

\title{Conservative and Stable Degree Preserving SBP Operators for Non-Conforming Meshes}

\author{Lucas Friedrich\thanks{Mathematical Institute, University of Cologne, Cologne, Germany}\and David C.~Del Rey Fern\'andez\thanks{Institute for Aerospace Studies, University of Toronto, Toronto, Canada}\and Andrew R. Winters\footnotemark[1]\and Gregor J. Gassner\footnotemark[1]\and David W. Zingg\footnotemark[2]\and Jason Hicken\thanks{Department of Mechanical, Aerospace, and Nuclear Engineering, Rensselaer Polytechnic Institute, Troy, New York, USA}}
\date{}
\maketitle
%

\begin{abstract}
Non-conforming numerical approximations offer increased flexibility for applications that require high resolution in a localized area of the computational domain or near complex geometries. Two key properties for non-conforming methods to be applicable to real
world applications are conservation and energy stability. The summation-by-parts (SBP) property, which
certain finite-difference and discontinuous Galerkin methods have, finds success for the
numerical approximation of hyperbolic conservation laws, because the proofs of
energy stability and conservation can discretely mimic the continuous analysis of partial
differential equations. In addition, SBP methods can be developed with high-order accuracy, which is
useful for simulations that contain multiple spatial and temporal scales. However,
existing non-conforming SBP schemes result in a reduction of the overall degree of the scheme, which leads to a reduction in the order of the solution error. This
loss of degree is due to the particular interface coupling through a
simultaneous-approximation-term (SAT). We present in this work a novel class of SBP-SAT operators
that maintain conservation, energy stability, and have no loss of the degree of the scheme for non-conforming
approximations. The new \emph{degree preserving} discretizations require an ansatz that the norm matrix of the SBP
operator is of a degree $\geq 2p$, in contrast to, for example, existing finite difference SBP operators,
where the norm matrix is $2p-1$ accurate.
We demonstrate the fundamental properties of the new scheme with rigorous mathematical analysis as well as numerical verification.
\end{abstract}



\parskip 2 mm

\section{Introduction}

As we move closer to the advent of exascale high performance computing, it becomes increasingly evident
that the numerical simulation of real-world applications on such machines requires flexible and robust methods.
One approach that provides numerical efficiency and robustness is the combination of summation-by-parts (SBP) operators~\cite{Fernandez2014,Svard2014,DCDRF2014,Gassner2013,Carpenter1996}
with simultaneous-approximation-terms (SATs)~\cite{Carpenter1994,Carpenter1999,Nordstrom1999,Nordstrom2001b} for the weak imposition of boundary conditions and interface coupling.
These nodal based SBP-SAT schemes are advantageous as they are conservative, high-order, linearly~\cite{Fernandez2014,Svard2014} and nonlinearly stable~\cite{Gassner2016,Matteo2015,Carpenter2014,Fisher2013b} and are applicable to structured multiblock~\cite{Fernandez2014,Svard2014} or unstructured
meshes~\cite{Fernandez2016,HickenAIAA2016,Hicken2016,Hicken2015b}.

Many real world problems contain a wide range of length scales, and the efficient approximation of
such problems necessitates the ability to judiciously distribute degrees of freedom.
The goal is to construct discretizations that formally retain their convergence order across non-conforming interfaces and therefore, this paper is a first step in the construction of adaptive SBP-SAT schemes with such properties.

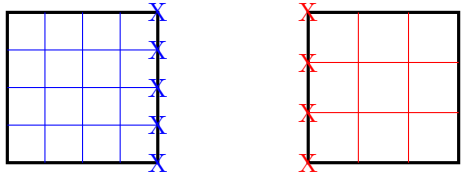
\begin{figure}[!ht]
\begin{center}
\begin{tikzpicture}[scale=0.5]
\draw[very thick] (-6,-2) -- (-2,-2) -- (-2,2) -- (-6,2) -- cycle;
\draw[very thick] (2,-2) -- (6,-2) -- (6,2) -- (2,2) -- cycle;
\color[rgb]{0,0,1}
\draw[] (-5,2) -- (-5,-2) ;
\draw[] (-4,2) -- (-4,-2) ;
\draw[] (-3,2) -- (-3,-2) ;
\draw[] (-6,1) -- (-2,1) ;
\draw[] (-6,0) -- (-2,0) ;
\draw[] (-6,-1) -- (-2,-1) ;
\node at (-2,-2) {X};
\node at (-2,-1) {X};
\node at (-2,0) {X};
\node at (-2,1) {X};
\node at (-2,2) {X};
\color[rgb]{1,0,0}
\draw[] (3.33333,2) -- (3.33333,-2) ;
\draw[] (4.66666,2) -- (4.66666,-2) ;
\draw[] (2,0.66666) -- (6,0.66666) ;
\draw[] (2,-0.66666) -- (6,-0.66666) ;
\node at (2,-2) {X};
\node at (2,-0.66666) {X};
\node at (2,0.66666) {X};
\node at (2,2) {X};
\end{tikzpicture}
		\caption{Two conforming elements with non-conforming nodal distributions along their common interface.}
		\label{fig:non}
\end{center}
\end{figure}

We concentrate on the development of non-conforming SBP-SAT finite difference methods due to their ability to arbitrarily assign degrees of freedom within elements. The construction of the SBP-SAT operators herein assumes that there is no subdivision of elements. That is, we focus on the interface coupling procedures where the elements are conforming but the distribution of nodes inside of the elements can vary, see Fig.~\ref{fig:non}. From here on the phrase \emph{non-conforming elements} refers to this type of non-conformity.

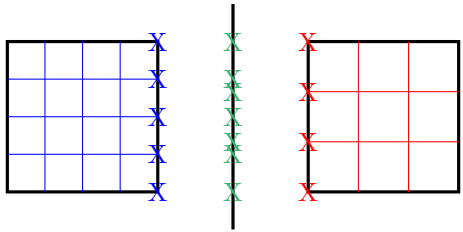
\begin{figure}[!ht]
\begin{center}
\begin{tikzpicture}[scale=0.5]
\draw[very thick] (-6,-2) -- (-2,-2) -- (-2,2) -- (-6,2) -- cycle;
\draw[very thick] (2,-2) -- (6,-2) -- (6,2) -- (2,2) -- cycle;
\draw[very thick] (0,3) -- (0,-3) ;
\color[rgb]{0,0,1}
\draw[] (-5,2) -- (-5,-2) ;
\draw[] (-4,2) -- (-4,-2) ;
\draw[] (-3,2) -- (-3,-2) ;
\draw[] (-6,1) -- (-2,1) ;
\draw[] (-6,0) -- (-2,0) ;
\draw[] (-6,-1) -- (-2,-1) ;
\node at (-2,-2) {X};
\node at (-2,-1) {X};
\node at (-2,0) {X};
\node at (-2,1) {X};
\node at (-2,2) {X};
\color[rgb]{1,0,0}
\draw[] (3.33333,2) -- (3.33333,-2) ;
\draw[] (4.66666,2) -- (4.66666,-2) ;
\draw[] (2,0.66666) -- (6,0.66666) ;
\draw[] (2,-0.66666) -- (6,-0.66666) ;
\node at (2,-2) {X};
\node at (2,-0.66666) {X};
\node at (2,0.66666) {X};
\node at (2,2) {X};
\color[rgb]{0.24,0.7,0.44}
\node at (0,-2) {X};
\node at (0,-1) {X};
\node at (0,-0.66666) {X};
\node at (0,0) {X};
\node at (0,0.66666) {X};
\node at (0,1) {X};
\node at (0,2) {X};
\end{tikzpicture}
		\caption{Two conforming elements with non-conforming nodal distributions and an intermediate grid.}
		\label{fig:nongrid}
\end{center}
\end{figure}

\textblue{Mattsson and Carpenter~\cite{Mattsson2010b} developed interpolation operators that result in conservative and stable schemes when incorporated
into the SBP-SAT framework. Here, the interpolation operators depend on the neighbouring non-conforming  element. Building on these ideas,
Kozdon and Wilcox~\cite{Kozdon2016} developed an approach that couples non-conforming
elements by first projecting the solution in adjacent elements onto an intermediate grid, see Fig.~\ref{fig:nongrid}, and then projecting back to the surface of the corresponding element. The main advantage of this approach is that one can construct interpolation/projection operators independently of the neighbouring elements as all element interfaces are projected to the same set of nodes. Throughout this paper we will adapt this idea and consider different types of intermediate grids.}

For clarity, we outline the four fundamental properties that we wish the new non-conforming SBP-SAT discretization to satisfy:
\begin{enumerate}
\item \emph{Conservation}: The numerical approximation must discretely recover the first moment, i.e. the integrals of the variables of the conservation law. For example, this ensures no loss of the total mass in time.\\[0.1cm]

\item \emph{Energy stability}: The numerical solution remains bounded by the initial conditions and boundary conditions of the problem. For a simple example of the continuous energy analysis consider the one-dimensional linear advection equation with unit wave
speed,
\begin{equation}\label{advection}
\frac{\partial \fnc{U}}{\partial t}=-\frac{\partial \fnc{U}}{\partial \xi},\quad \xi\in[0,1],\quad t\ge0,
\end{equation}
together with an initial condition and a Dirichlet boundary condition
\begin{equation}\label{advectionICBC}
\fnc{U}(\xi,0)=\fnc{U}_0(\xi),\quad \fnc{U}(\xil,t)=\fnc{G}_{L}(t).
\end{equation}
To prove that the problem defined by (\ref{advection}) and (\ref{advectionICBC}) is stable, we employ the energy method. First, we
multiply (\ref{advection}) by the solution and integrate over the domain to find
\begin{equation}\label{advection:E1}
\int\limits_{0}^{1}\fnc{U}\frac{\partial \fnc{U}}{\partial t}\mr{d}\xi=-\int\limits_{0}^{1} \fnc{U} \frac{\partial \fnc{U}}{\partial \xi}\,\mr{d}\xi.
\end{equation}
For the left-hand term of (\ref{advection:E1}), we bring $\fnc{U}$ into the temporal derivative and apply Leibniz's rule. For the right-hand term of (\ref{advection:E1}), we use integration-by-parts, resulting in
\begin{equation}
\frac{\mr{d}(\fnc{U},\fnc{U})}{\mr{d}t}=-\left.\fnc{U}^{2}\right|_{0}^{1},
\end{equation}
where we use the $\Ltwo$ inner product
\begin{equation}
(\fnc{U},\fnc{U}) \coloneqq  \int\limits_{0}^{1}\fnc{U}^2\,\mr{d}\xi.
\end{equation}
Finally, we apply the boundary condition, initial condition, integrate to the final time $t=T$
and rearrange to obtain
\begin{equation}\label{eq:energy}
(\fnc{U},\fnc{U}) \leq (\fnc{U}_0,\fnc{U}_0) + \int_{0}^{T}\fnc{G}_{L}^{2}\,\mr{d}t,
\end{equation}
which shows that the solution is bounded by the initial condition and the boundary
data. Therefore, the problem is stable in the sense of Hadamard \cite{Gustafsson1996}. That is, the solution is continuously
dependent on the initial data of the problem.\\[0.1cm]

\item \emph{High-order}: The scheme can differentiate polynomials with high degree exactly; this implies that the order of the solution error is high-order.\\[0.1cm]

\item \emph{Degree preservation}: In this work degree refers to the highest degree, say $p$, of the monomial for which the differentiation matrix is exact. The constituent components of the SBP scheme, i.e., the derivative operators and SATs, remain of degree $p$ when applied to non-conforming meshes. 
\end{enumerate}

The first three properties are immediately available for non-conforming existing SBP finite difference schemes from the work of Kozdon and Wilcox \cite{Kozdon2016} and Mattsson and Carpenter \cite{Mattsson2010b}. Unfortunately, it was shown by Lundquist and Nordstr\"om~\cite{Nordstrom2015} that existing degree $p$ SBP schemes, when applied to non-conforming meshes, result in SATs that lose one degree of accuracy compared to the conforming case. Therefore, existing SBP schemes are not \emph{degree preserving} and thus violate the fourth desired.

The primary objective of this work is to construct SBP operators such that the resulting scheme is degree preserving. It turns out that the essential idea is to consider SBP operators with norm matrices that are of at least degree $2p$, which allows the construction of degree preserving SATs. The second objective is to present a generalized construction of SATs that leads to degree preserving schemes. 

\subsection{Notational conventions}\label{sec:notational}

To define precisely the mathematical analysis of the new SBP schemes we first introduce the notation that we use throughout this work, adapted from Hicken et al.~\cite{Hicken2016}. To reduce the notational complexity, the presentation is restricted to two-dimensional operators.
Vectors are denoted with lower case bold letters, for example $\bm{\xi}=\left[\xi_{1},\dots,\xi_{N}\right]\Tr$, while matrices are presented using
capital letters with sans-serif font, for example, \mat{H}. Script letters are used to denote continuous functions on a
specified domain, e.g.,
\begin{equation}\mathcal{U}(\xi,\eta,t)\in \Ltwo\left(\hat{\Omega}\times [0,T]\right),\end{equation}
denotes a square integral function on the spatial domain $\hat{\Omega} = \left[\xi_{1},\xi_{\Nxi}\right]\times[\eta_{1},\eta_{\Neta}]$. 
The restriction of functions onto a computational mesh is denoted with lower case bold font.
For example, the restriction of $\fnc{U}$ onto a grid of $N$ nodes in each spatial direction, $S_{\hat{\Gamma}}=\left\{\left(\xi_{i},\eta_{i}\right)\right\}_{i=1}^{N}$, is given by
\begin{equation}
\bm{u}\coloneqq\left[\fnc{U}\left(\xi_{1},\eta_{1}\right),\dots,\fnc{U}\left(\xi_{N},\eta_{N}\right)\right]\Tr.
\end{equation}
Monomial basis functions are used in the proofs of the underlying properties of the new non-conforming SBP-SAT schemes. 
We define the monomials for the tensor-product basis by
\begin{equation}
\Pk(\xi,\eta)\coloneqq  \xi^{i}\eta^{j},
\end{equation}
where $k$ is uniquely associated with a monomial with individual powers of
\begin{equation}\label{convk}
k=j(p+1)+i+1,\quad i,j=0,\ldots,p.
\end{equation}
The cardinality of a basis for the polynomials under consideration is given by
\begin{equation}
\nmin{p}\coloneqq  (p+1)^d.
\end{equation}
where $d$ is the spatial dimension.

%

The monomials and their derivatives evaluated at the computational nodes are represented by
\begin{equation}\label{eq:monomial}
\bpk\coloneqq \left[\Pk(\xi_{1},\eta_{1}),\dots,\Pk(\xi_{N_\xi},\eta_{N_\eta})\right]\Tr,
\end{equation}
and
\begin{equation}
\dxibpk\coloneqq \left[\frac{\partial \Pk}{\partial \xi}(\xi_{1},\eta_{1}),\dots,\frac{\partial \Pk}{\partial \xi}(\xi_{N_\xi},\eta_{N_\eta})\right]\Tr.
\end{equation}
Typically, SBP methods are built from tensor-product operators on Cartesian grids for which the nodal distribution is
defined by two vectors
\begin{equation}
\bm{\xi}=\left[\xi_{1},\dots,\xi_{\Nxi}\right]\Tr,\quad
\bm{\eta}=\left[\eta_{1},\dots,\eta_{\Neta}\right]\Tr.
\end{equation}
Thus the total number of nodes is $N=\Nxi\Neta$.

Having defined $\bm{\xi}$ and $\bm{\eta}$, the projection of the monomials in \eqref{eq:monomial} onto the grid
is constructed by
\begin{equation}
\bpk = \bm{\xi}^{i}\otimes\bm{\eta}^{j},
\end{equation}
where $\otimes$ is the tensor-product operator, and $i$ and $j$ can be found from $k$
using (\ref{convk}).
Furthermore, for any vector $\bm{v}$ we define
\begin{equation}
\bm{v}^{r} \coloneqq  \left[v_{1}^{r},\dots,v_{N}^{r}\right]\Tr,
\end{equation}
where if $r<0$ then $\bm{v}^{r}$ is a vector of zeros.

The degree of a SBP operator is the degree of the monomial
for which the differentiation matrix operator is exact at all nodes. For example, assume a set of computational nodes $\bm{\xi}$ and a one-dimensional differentiation operator $\Dxione$. If
\begin{equation}\label{eg: Definition Degree}
\Dxione\bm{\xi}^k=k\bm{\xi}^{k-1},\qquad\text{ for }k=0,\dots,p,
\end{equation}
then the degree of the operator is $p$.

\subsection{Review of summation-by-parts operators}\label{sec:SBP}

Next, we give a brief introduction to the methodology of SBP methods. In essence, the goal of SBP operators is to discretely mimic the integration-by-parts property of first and higher derivatives~\cite{Fernandez2015,Mattsson2014,Mattsson2012,Mattsson2004b}. We demonstrate, in a one-to-one fashion, how the steps of the continuous energy analysis are mimicked at the semi-discrete level with SBP operators to prove stability of semi-discrete approximations. Complete details on SBP operators and the energy method can be found in, e.g., \cite{Fernandez2014,Svard2014}.

A numerical scheme is said to be stable if the semi-discrete solution is bounded in terms of the initial and boundary data of the problem. To demonstrate the discrete energy stability we define a first-derivative SBP operator that is applicable to general one-dimensional nodal distributions~\cite{DCDRF2014}.
\begin{definition}\label{DEFSBP}
{\bf One-dimensional summation-by-parts operator for the first derivative:} A matrix operator
$\Dxione \in\mathbb{R}^{N\times N}$ is an SBP operator of degree $p$ which approximates the first derivative $\frac{\partial}{\partial \xi}$
on the nodal distribution $\bm{\xi}$ in computational space, if
\begin{enumerate}
\item $\Dxione\xik{k}= \left(\Mxione\right)^{-1}\Qxione\xik{k}=\left(\Mxione\right)^{-1}\left(\Sxione+\frac{1}{2}\ExioneD\right)\xik{k}= k\xik{k-1}$, $k=0,1,\dots,p$.\\[0.1cm]
\item $\Mxione$ denotes the norm matrix and is symmetric positive definite.\\[0.1cm]
\item $\ExioneD$ is symmetric and  $\Sxione$ is skew-symmetric. Thus, $\Qxione+\left(\Qxione\right)\Tr=\ExioneD$.\\[0.1cm]
\item $\left(\xik{i}\right)\Tr\ExioneD\xik{j}=\xi_{\Nxi}^{i+j}-\xi_{1}^{i+j},\quad i,j=0,1,\dots,r, r\ge p$.\\[0.1cm]
\end{enumerate}
\end{definition}
\begin{rmk}
In this paper we consider only SBP operators where the norm matrix $\Mxione$ is diagonal.
\end{rmk}
To impose boundary conditions or inter-element coupling we use the concept of SATs. The definition of the SAT requires the decompositions of $\ExioneD$ into contributions associated with the boundaries of each element. In this paper, we consider SBP operators that contain nodes at each of the boundaries.
A simple decomposition of $\ExioneD$ is given by \cite{DCDRF2014}
\begin{equation}\label{decompEx}
\ExioneD=\txir\txir\Tr-\txil\txil\Tr,
\end{equation}
where
\begin{equation}\label{condtxLtxR}
\txir=\left[0,\dots,0,1\right]\Tr,\quad\txil=\left[1,0,\dots,0\right]\Tr.
\end{equation}

The discretization of the spatial terms in the one-dimensional linear advection equation (\ref{advection}) using SBP operators results in
\begin{equation}\label{eq:disc}
\frac{\mr{d}\bm{u}}{\mr{d}t} = -\Dxione\bm{u}+\sigma\left(\Mxione\right)^{-1}\left(\txil\txil\Tr\bm{u}-\fnc{G}_{L}\txil\right).
\end{equation}
The last term in \eqref{eq:disc} is the SAT,
\cite{Carpenter1994,Carpenter1999,Nordstrom1999,Nordstrom2001b}, which weakly imposes
the boundary condition \eqref{advectionICBC}. Note that the SAT incorporates the yet-to-be-determined scalar stability parameter $\sigma$.

Next, we show how the discrete energy method constrains the choice of
$\sigma$ such that the resulting discretization is stable. The symmetric positive definite norm matrix $\Mxione$
generates a discrete $\Ltwo$ inner product \cite{DCDRF2014,Hicken2013} such that
\begin{equation}\label{eq:discreteL2}
\langle \bm{u},\bm{v}\rangle \coloneqq  \bm{v}\Tr\Mxione\bm{u}\approx\int\limits_{0}^{1}\fnc{U}\fnc{V}\,\mr{d}\xi = (\fnc{U},\fnc{V}).
\end{equation}
Multiplying \eqref{eq:disc} by $\bm{u}\Tr\Mxione$ from the left, using the third property of the SBP operator, and rearranging the equation, we obtain
\begin{equation}\label{advectionE4}
\frac{\mr{d}\langle\bm{u},\bm{u}\rangle}{\mr{d}t}=-\bm{u}\Tr\left(\txir\txir\Tr-\txil\txil\Tr\right)\bm{u}+2\sigma\bm{u}\Tr\left(\txil\txil\Tr\bm{u}-\fnc{G}_{L}\txil\right).
\end{equation}
Finally, completing the square on the right-hand side,
integrating up to a final time $T$, and including the initial condition, (\ref{advectionE4}) becomes
\begin{equation}
\langle\bm{u},\bm{u}\rangle\leq \langle\bm{u}_0,\bm{u}_0\rangle-\gamma\Gamma^2\int_{0}^{T}\fnc{G}_{L}^{2}\,\mr{d}t,
\end{equation}
where we introduce two additional constants $\Gamma = \frac{-\sigma}{2\sigma+1}$ and $\gamma=2\sigma+1$. For the scheme to remain stable we require $\gamma<0$. We know that $\sigma =-1$ for the scheme to remain conservative \cite{Fernandez2014}. With this selection of $\sigma$ we exactly mimic the continuous energy result \eqref{eq:energy} term-by-term and obtain the discrete energy
\begin{equation}
\langle\bm{u},\bm{u}\rangle\leq \langle\bm{u}_0,\bm{u}_0\rangle+\int_{0}^{T}\fnc{G}_{L}^{2}\,\mr{d}t.
\end{equation}
The extension of one-dimensional SBP operators to
two spatial dimensions is accomplished with tensor products. 
%
\begin{equation}\label{eq:Dtensor}
\Dxi = \Dxione\otimes\Ieta,\quad\Deta=\Ixi\otimes\Detaone,
\end{equation}
where $\Ixi$ and $\Ieta$ denote the identity matrix of size $\Nxi$ and $\Neta$, respectively.

The two-dimensional SBP operators satisfy
\begin{equation}\label{eq:2DOperators}
\Dxi\coloneqq \M^{-1}\Qxi=\M^{-1}\left(\Sxi+\frac{1}{2}\Exi\right),\quad\Deta\coloneqq \M^{-1}\Qeta=\M^{-1}\left(\Seta+\frac{1}{2}\Eeta\right),
\end{equation}
where
\begin{equation}\label{set:Operators}
\begin{split}
&\M\coloneqq \Mxione\otimes\Metaone,\\[0.1cm]
 &\Qxi\coloneqq \Qxione\otimes\Metaone,\quad\Qeta\coloneqq \Mxione\otimes\Qetaone,\quad\Sxi\coloneqq \Sxione\otimes\Metaone,\quad\Seta\coloneqq \Mxione\otimes\Setaone,\\[0.1cm]
 &\Exi\coloneqq \ExiN-\Exione,\quad\ExiN\coloneqq \RxiN\Tr\Metaone\RxiN,\quad\Exione\coloneqq \Rxione\Tr\Metaone\Rxione,\\[0.1cm]
 &\Eeta\coloneqq \EetaN-\Eetaone,\quad\EetaN\coloneqq \RetaN\Tr\Mxione\RetaN,\quad\Eetaone\coloneqq \Retaone\Tr\Mxione\Retaone,\\[0.1cm]
 &\RxiN \coloneqq  \txiN\Tr\otimes\Ieta,\quad\Rxione \coloneqq  \txione\Tr\otimes\Ieta,\quad
\RetaN \coloneqq  \Ixi\otimes\tetaN\Tr,\quad\Retaone \coloneqq  \Ixi\otimes\tetaone\Tr.
\end{split}
\end{equation}
%
%
It is useful to note that the constituent matrices of \eqref{set:Operators} are approximations of the following bilinear forms~\cite{Hicken2016}:
\begin{equation}\label{eq:2.5}
\begin{split}
&\bm{v}\Tr\M\bm{u}\approx\int_{\Omega}\fnc{V}\fnc{U}\mr{d}\Omega,
\quad\bm{v}\Tr\Qxi\bm{u}\approx\int_{\Omega}\fnc{V}\frac{\partial\fnc{U}}{\partial \xi}\mr{d}\Omega,
\quad\bm{v}\Tr\Qeta\bm{u}\approx\int_{\Omega}\fnc{V}\frac{\partial\fnc{U}}{\partial \eta}\mr{d}\Omega,\\
&\bm{v}\Tr\Exi\bm{u}\approx\oint_{\Gamma}\fnc{V}\fnc{U}\nxi\mr{d}\Gamma,\quad\bm{v}\Tr\Eeta\bm{u}\approx\oint_{\Gamma}\fnc{V}\fnc{U}\neta\mr{d}\Gamma.
\end{split}
\end{equation}
In particular, we see that $\Exi$ and $\Eeta$ are bilinear forms approximating surface integrals. Later we use this
interpretation to clarify the meaning of the SATs.

The remainder of the paper is organized as follows: Section \ref{sec:Construction} builds a non-conforming discretization with two-dimensional SBP operators. Next, we create a novel set of degree preserving SBP operators in Sec. \ref{sec:examples}. These new operators create a non-conforming numerical scheme that is conservative, energy stable, high-order and degree preserving, as shown in Sec. \ref{sec:proof}. Numerical results are presented in Sec. \ref{sec:results} to support the theoretical findings. Our concluding remarks are given in the final section.

%
%


\section{Construction of coupling SATs for non-conforming elements}\label{sec:Construction}
Mattsson and Carpenter \cite{Mattsson2010b} developed interpolation operators to project the approximate solution of a finite difference scheme from one element to another in a stable way. In \cite{Kozdon2016}, Kozdon and Wilcox expanded this conservative and stable finite difference scheme to non-conforming approximations where each element can have an arbitrary number of computational nodes with the introduction of intermediate \emph{glue grids}. 
However, the non-conforming scheme of Kozdon and Wilcox does not preserve the degree of the conforming approximation \cite{Kozdon2016}; indeed, Lundquist and Nordstr\"om \cite{Nordstrom2015} proved that degree preserving SBP schemes cannot be constructed using existing SBP operators --- we prove this in Thm \ref{thrm:tpm1} and demonstrate this numerically in Sec. \ref{subsec:Coupling Classic}.

The essential idea of the proposed non-conforming scheme is to project the solution onto a new set of nodes at an interface,
and construct the SATs on these intermediate nodes.
%

For the analysis, we focus on the linear two-dimensional constant-coefficient advection equation
\begin{equation}\label{LC}
\frac{\partial\fnc{U}}{\partial t}+\betax\frac{\partial\fnc{U}}{\partial x}
+\betay\frac{\partial\fnc{U}}{\partial y} = 0,\quad(x,y)\in\Omega,\quad t\ge 0,
\end{equation}
where $\betax$ and $\betay$ are the constant wave speeds. The hyperbolic equation \eqref{LC} is subject to initial and boundary conditions
\begin{equation}\label{ICBC}
\fnc{U}(x,y,0) = \fnc{U}_0(x,y),\quad\fnc{U}(x,y,t) = \fnc{G}(x,y,t),\quad\forall(x,y)\in\Gamman.
\end{equation}
The boundary, $\Gamma$, of the domain
$\Omega$ is decomposed into two surfaces
$\Gamman=\left\{(x,y)\in\Gamma\,|\,(\betax\nx+\betay\ny<0)\right\}$,
$\Gammap=\Gamma/\Gamman$, and $\nx$, $\ny$, are the $x$ and $y$ components of
the outward pointing unit normal to $\Gamma$.

In order to use tensor-product SBP operators, it is first necessary to partition
the domain $\Omega$ into $K$ non-overlapping quadrilateral elements $\Omega_{i}$, such that
\begin{equation}
\Omega = \bigcup\limits_{i = 1}^{K}\Omega_{i}.
\end{equation}
The tensor-product SBP operators are defined on a regular Cartesian computational grid; therefore, on each element, the PDE (\ref{LC}) is mapped from physical coordinates, $(x,y)\in\Omega_{i}$
to computational coordinates $(\xi,\eta)\in\Ohat$, resulting in
\begin{equation}\label{CLCS}
\frac{\partial\left(\fnc{J}\fnc{U}\right)}{\partial t}+\frac{\partial\left(\lamxi\fnc{U}\right)}{\partial \xi}
+\frac{\partial\left(\lameta\fnc{U}\right)}{\partial \eta}=0,
\end{equation}
where
\begin{equation}\label{lams}
\lamxi = \betax\frac{\partial y}{\partial \eta}-\betay\frac{\partial x}{\partial \eta},\quad
\lameta = -\betax\frac{\partial y}{\partial \xi}+\betay\frac{\partial x}{\partial \xi},\quad
\fnc{J}= \frac{\partial x}{\partial\xi}\frac{\partial y}{\partial\eta}
-\frac{\partial x}{\partial \eta}\frac{\partial y}{\partial \xi}.
\end{equation}

\begin{rmk}
For simplicity, we restrict the discussion in this paper to only consider affine maps. Therefore the mapping terms, for
example $\fnc{J}$, are constants. The extension to curvilinear elements is straightforward and is nicely outlined in Kozdon and Wilcox \cite{Kozdon2016}.
\end{rmk}

\subsection{Conforming discretization}\label{sec:Conforming Discretization}
Before presenting SATs that lead to stable and conservative schemes for non-conforming elements,
it is instructive to examine SATs for the conforming case. This will clarify the meaning of the
SATs. Also, it allows us to reformulate the SATs using numerical fluxes as is commonly done in discontinuous
Galerkin methods, which is the form we will use throughout the paper.

\textred{We consider the coupling between a left element, $L$, and a right element, $R$, that are conforming
and share a vertical interface. The discretization in element $L$, considering only the coupling
at the shared interface, is given as
\begin{equation}\label{eq:discLcon}
\begin{split}
\frac{\mr{d}\uA}{\mr{d}t}+\J\lamxi\Dxi\uA+\J\lameta\Deta\uA =\J\M^{-1}\left(\lamxi\ExiN\uA-\fs\right),
\end{split}
\end{equation}
where $\bm f^{*}$ denotes the numerical flux function. The term on the right-hand side of (\ref{eq:discLcon}) is the SAT that couples the solution in element $L$ to
the solution in element $R$.
Similarly, the discretization in element $R$ is given as
\begin{equation}\label{eq:discRcon}
\begin{split}
\frac{\mr{d}\uB}{\mr{d}t}+\J\lamxi\Dxi\uB+\J\lameta\Deta\uB =-\J\M^{-1}\left(\lamxi\Exione\uB-\fs\right).
\end{split}
\end{equation}
where, for simplicity, we assume that both elements have identical nodal distributions.
We define the numerical flux functions
\begin{equation}\label{eq:centralflux}
\fs \coloneqq  \frac{1}{2}\left(\lamxi\ExiN\uA+\lamxi\Exione\uB\right)-\frac{\sigma|\lamxi|}{2}\left(\Exione\uR-\ExiN\uL \right).\\
\end{equation}
If we set $\sigma=0$, then we recover a central numerical flux. For $\sigma=1$ we obtain an upwind SAT. Note that the definition of the numerical flux includes the surface mass matrix and therefore its definition differs from that in \cite{Gassner2013}, where the numerical flux is defined by
\begin{equation}
\fsGas \coloneqq  \frac{1}{2}\left(\lamxi\RxiN\uA+\lamxi\Rxione\uB\right)-\frac{\sigma|\lamxi|}{2}\left(\Rxione\uR-\RxiN\uL \right).
\end{equation}}
To clarify the role of the SATs in the discretization, we recast (\ref{eq:discLcon}) and (\ref{eq:discRcon}) into weak
form. \textred{To do so, we consider a computational grid with an appropriate
functional space $\fnc{V}\in\mathscr{V}$}. Then we pre-multiply (\ref{eq:discLcon}) by
$\bm{v}_{L}\Tr\M$ and (\ref{eq:discRcon}) by $\bm{v}_{R}\Tr\M$, where $\bm{v}_{L}$ and $\bm{v}_{R}$
are basis functions of $\fnc{V}$ on to the nodes of the elements $L$ and $R$, respectively. Thus the
problem defined by (\ref{eq:discLcon}) becomes: Find $\uA$ such
that for all $\bm{v}_L\in\mathbb{R}^{N}$
\begin{equation}\label{eq:discLconfw}
\begin{split}
&\bm{v}_{L}\Tr\M\frac{\mr{d}\uA}{\mr{d}t}+\J\lamxi\bm{v}_{L}\Tr\Qxi\uA+\J\lameta\bm{v}_{L}\Tr\Qeta\uA =\J\bm{v}_{L}\Tr\left(\lamxi\ExiN\uA-\fs\right),
\end{split}
\end{equation}
where a similar weak problem holds for the solution in the right element $\uB$. From Section \ref{sec:SBP} we know that each matrix of an
SBP operator is an approximation to a certain bilinear form. Therefore, we see that
\begin{equation}\label{eq:discLconfwint}
\begin{split}
&\overbrace{\bm{v}_{L}\Tr\M\frac{\mr{d}\uA}{\mr{d}t}}^{\approx\int_{\Ohat_{L}}\fnc{V}\frac{\partial\fnc{U}}{\partial t}\mr{d}\Ohat}+
\overbrace{\J\lamxi\bm{v}_{L}\Tr\Qxi\uA+\J\lameta\bm{v}_{L}\Tr\Qeta\uA}^{\approx\J\int_{\Ohat_{L}}\fnc{V}\left(\lamxi\frac{\partial\fnc{U}}{\partial\xi}+\lameta\frac{\partial\fnc{U}}{\partial\eta}\right)\mr{d}\Ohat}
=\overbrace{\J\bm{v}_{L}\Tr\left(\lamxi\ExiN\uA-\fs\right)}^{\approx\J\oint_{\Ghat}\fnc{V}\left(\lamxi\fnc{U}-\fnc{F}_{Gas}^{*}\right)\nxiL\mr{d}\Ghat},
\end{split}
\end{equation}
where $\nxiL$ is the $\xi$ component of the outward pointing unit normal at the shared interface for element $L$. Reformulating the
SBP-SAT discretization in weak form\footnote{Here we have not transferred the action of the derivative onto the test function; in the present context these two forms are algebraically equivalent as a result of constant grid metrics \cite{Kopriva2010}.} clarifies that the components of the SAT are surface integrals. The $\mat{E}$ matrices in \eqref{eq:discLcon} and \eqref{eq:discRcon} project the solution onto the nodes of the common interface and approximate the surface integral. The construction of SATs for non-conforming elements is related to this idea.

\subsection{Non-conforming discretization}\label{sec:Non-Conforming Discretization}
Now we consider the case of two elements that do not have conforming nodal distributions.
For coupling the solution between elements we introduce an intermediate grid. \textred{On this intermediate grid we define
a set of nodes $\{(\xi_{\Gamma,i},\eta_{\Gamma,i})\}_{i=1}^{N_\Gamma}$
with a corresponding symmetric positive definite surface-norm matrix $\MGam$ of sufficient accuracy}. In contrast to \cite{Kozdon2016,Mattsson2010b}, we do not project back to each respective element. Rather we directly use the node distribution of the intermediate grid to construct approximations of the surface integrals.

\textred{From the conforming discretizations \eqref{eq:discLcon} and \eqref{eq:discRcon}, we generalize to a non-conforming discretization
\begin{equation}\label{eq:discA}
\begin{split}
&\frac{\mr{d}\uA}{\mr{d}t}+\JA\lamxiA\DxiA\uA+\JA\lametaA\DetaA\uA =\frac{1}{\fnc{J}_L}\HA^{-1}\left(\lamxiA\ExiNL\uA-\fsL\right),
\end{split}
\end{equation}
in the left element and
\begin{equation}\label{eq:discB}
\begin{split}
&\frac{\mr{d}\uB}{\mr{d}t}+\JB\lamxiB\DxiB\uB+\JB\lametaB\DetaB\uB =-\frac{1}{\fnc{J}_R}\HB^{-1}\left(\lamxiB\ExioneR\uB-\fsR\right),
\end{split}
\end{equation}
in the right element. Note that although the nodal distributions differ on each element at the
interface, $\lamxiL=\lamxiR$ and $\lametaL=\lametaR$ because the elements are conforming (but not the nodal distribution) and we consider an affine mapping. For non-conforming discretizations the numerical flux is defined as
\begin{equation}\label{eq:ustarLR}
\begin{split}
&\fsL=\frac{1}{2}\left(\lamxi\ExiNL\uL+\lamxi\RL\Tr\MGam\RR\uR \right)-\frac{\sigma|\lamxi|}{2}\left(\RL\Tr\MGam\RR\uR-\RL\Tr\MGam\RL\uL \right),\\
&\fsR=\frac{1}{2}\left(\lamxi\ExioneR\uR+\lamxi\RR\Tr\MGam\RL\uL \right)-\frac{\sigma|\lamxi|}{2}\left(\RR\Tr\MGam\RR\uR-\RR\Tr\MGam\RL\uL \right),
\end{split}
\end{equation}
with $\lamxi\coloneqq \lamxiL=\lamxiR$ and $\lameta\coloneqq \lameta=\lametaR$. Note that the coupling procedure in the non-conforming approximation is more complicated. We now include the norm matrix from the intermediate grid $\MGam$ and introduce the new \emph{projection operators} $\RL$ and $\RR$. In this paper we consider projection operators that are constructed using tensor products. Therefore we can rewrite the projection operators in terms of
\begin{equation}
\RL=\exiNL\Tr\otimes\RLone\quad\RR=\exioneR\Tr\otimes\RRone.
\end{equation}
Here the projection operators $\RLone,\RRone$ project the solution from the element interface on to the intermediate grid. Throughout this section, all derivations are in one-dimension as we focus on interfaces.
\begin{rmk}
Considering tensor products the proposed SATs are equivalent to the approach of Kozdon and Wilcox \cite{Kozdon2016}, as shown in Appendix \ref{sec:KW}. We concentrate on tensor-product projection operators in order to demonstrate that the new degree preserving SBP operators can be equally applied in the context of other coupling procedures, e.g., \cite{Kozdon2016,Mattsson2010b}.
\end{rmk}
In order to construct the accuracy conditions for the SAT, it is necessary to map the monomials in one computational space to the other. We therefore need to introduce a set of nodes $\etaL$ and $\etaR$, which denote the nodal distribution of the left and right element in one dimension along the interface. Furthermore, we define $\etaGam$ which denotes the corresponding nodes of the intermediate grid.
In order to construct a tensor-product projection operator of degree $p$, the operators must satisfy
\begin{equation}\label{eq:RLRR}
\begin{split}
\RLone\etaLk&=\etaGamk,\\
\RRone\etaRk&=\etaGamk,
\end{split}
\end{equation}
for $k=0,1,\dots, p$. For the SATs to be of the same degree as the derivative operators, the following relations must be
satisfied:
\begin{equation}\label{eq:do}
\begin{split}
\RLoneT\MGam\RRone\etaRk&=\MetaLone\etaLk,\\
\RRoneT\MGam\RLone\etaLk&=\MetaRone\etaRk,
\end{split}
\end{equation}
for $k=0,1,\dots,p$. These conditions ensure that the SAT is zero when its arguments are polynomials of degree less than or equal to $p$ in the $\xi$-direction.
Including \eqref{eq:RLRR} in \eqref{eq:do}, we arrive at
\begin{equation}\label{eq:dofinal}
\begin{split}
\RLoneT\MGam\etaGamk&=\MetaLone\etaLk,\\
\RRoneT\MGam\etaGamk&=\MetaRone\etaRk,
\end{split}
\end{equation}
for $k=0,1,\dots,p$.
We see that the conditions on $\RL$ are independent of those on $\RR$. Therefore, the
construction of the projection operators can be performed independently. A brief derivation of the projection operators is provided in Appendix \ref{sec:Projection}.
Now we focus on the corresponding SBP operators for the non-conforming problem. Here, we define the degree of a norm matrix. Consider two monomials $\eta^i$ and $\eta^j$ with $i,j=0,\dots,p$. Referring to the norm matrix, we understand the term \emph{degree of $\tilde p$} as exactness in the following sense
\begin{equation}\label{eg:degree norm}
\left(\left(\bm{\eta^{i}}\right)\Tr\Mxione\bm{\eta^{j}}\right)=\left(\int_{\eta_1}^{\eta_{\Neta}}\eta^{i+j}\mr{d}\eta\right),\quad\text{ for }\quad i+j\le \tilde p.
\end{equation}
Similarly, the interpretation of degree is applied to the remaining constituent matrices of the SBP operators, see \eqref{eq:2.5}.
Most tensor-product SBP operators are implicitly constructed such that the norm matrices are of degree $2p-1$~\cite{DCDRF2014,Hicken2013}. For example,
existing finite difference operators \cite{Hicken2013} or the
discontinuous Galerkin spectral element SBP operators constructed on the
Legendre-Gauss-Lobatto (LGL) nodes~\cite{DCDRF2014,Gassner2013} have this property. A major drawback of such operators is that it is
not possible to construct $\RLone$ and $\RRone$ as in \eqref{eq:RLRR} and \eqref{eq:dofinal} such that the resulting SATs retain the accuracy of
the derivative operators, as proven in the following Theorem, where we adapt the ideas of Lundquist and Nordstr\"om \cite{Nordstrom2015}.
\begin{thm}\label{thrm:tpm1}
Given a degree $p$ SBP operator $\DxioneL$ with norm matrix $\MxiLone$ of degree $2p-1$, assuming an intermediate grid with a norm matrix $\MGam$ of degree $\ge 2p-1$ such that the matrices $
\ML$ and $\MGam$ are different in terms of the norm matrices having different errors, then it is
not possible to construct a projection operator $\RL$ that is of degree $p$.
\end{thm}
\begin{pf}
We assume that it is possible to construct operators of degree $p$ and seek a contradiction. Therefore, we have the degree $p$ projection operator $\RL$ with the properties
\begin{equation}\label{eq:accP}
\begin{split}
\RLone\etaLk&=\etaGamk, \\
\RLoneT\MGam\etaGamk&=\MxiLone\etaLk,
\end{split}
\end{equation}
for $k=1,\dots,p$.
Consider the case $k=p$. The norm matrix $\MxiLone$ is of degree $2p-1$, therefore
\begin{equation}\label{eq:accH}
\left(\etaL^p\right)\Tr\MxiLone\etaL^p \neq \int_{\eta_1}^{\eta_{\Neta}}\eta^{2p}\mr{d}\eta.
\end{equation}
By re-deriving \eqref{eq:accP} we get
\begin{align}
\RLoneT\MGam\etaGamk&=\MxiLone\etaLk,\notag\\
\Leftrightarrow\left(\etaL^p\right)\Tr\RLoneT\MGam\etaGam^p&=\left(\etaL^p\right)\Tr\MxiLone\etaL^p,\notag\\
\stackrel{\eqref{eq:accP}}{\Leftrightarrow}\left(\etaGam^p\right)\Tr\MGam\etaGam^p&=\left(\etaL^p\right)\Tr\MxiLone\etaL^p.\label{eq:Nordstrom}
\end{align}
However, by \eqref{eq:accH} the norm matrix $\MxiLone$ cannot integrate $\eta^{2p}$ exactly. Also, if $\MGam$ is of degree $2p-1$, then the left-hand side of \eqref{eq:Nordstrom} also does not integrate exactly. The error terms for each quadrature rule are different, so the equality \eqref{eq:Nordstrom} cannot hold in general. For the case where the degree of $\MGam$ is larger than $2p-1$, then the equality \eqref{eq:Nordstrom} cannot hold since the left-hand side performs exact integration of $\eta^{2p}$, whereas the right hand side produces an error.  \qed
\end{pf}}
Theorem \ref{thrm:tpm1} holds the key to constructing SBP operators for which projection operators
of degree $p$ can be built which lead to stable schemes. Namely, if one builds SBP operators that
have norm matrices that are of degree $\ge 2p$ then the equality that failed in the proof of Thm. \ref{thrm:tpm1} can be
satisfied for degree $p$ monomials. This is not only related to the presented discretization in \eqref{eq:discA}, \eqref{eq:discB}, but also holds for \cite{Mattsson2010b,Kozdon2016}, as proven by Lundquist and Nordstr\"om \cite{Nordstrom2015}.
We now prove under what conditions SBP operators of degree $p$ with norm matrices of at least degree $2p$, which we denote {\it degree preserving SBP operators}, exist.
\begin{thm}
The existence of a quadrature rule with positive weights of degree $\ge 2p$ is necessary and sufficient for the existence of a degree preserving SBP operator of degree $p$.
\end{thm}
\begin{pf}
The proof follows identically from that given in Del Rey Fern\'andez et al. \cite{DCDRF2014} for Thm. $2$ and is omitted for brevity. \qed
\end{pf}
Next, we construct degree preserving SBP operators.

\section{Example of a degree preserving SBP operator}\label{sec:examples}
In this section, we explicitly construct degree preserving SBP operators $(\Mxione,\Dxione)$ in one dimension. The two-dimensional SBP operator can be derived from \eqref{eq:Dtensor}. The degree of $\Dxione$ is $p$ and the degree of $\Mxione$ is $2p$. These operators are derived on a master element $\hat\Omega:=[-1,1]$. Here, $\Nxi$ denotes the number of nodes.

First, assume we want to construct a classical finite difference SBP operator, denoted
\textit{classical FD-SBP} operator $(\Mclas,\Dclas)$ with
$p\le 4$. \textred{Here, classical FD-SBP operators refer to SBP operators where the choice of the number of nodes $\Nxi$ is arbitrary assuming that $\Nxi$ is larger than a minimum number of nodes, see \cite{Fernandez2014,Fernandez2015}. For such operators the norm matrices are of degree $2p-1$. To construct these operators}, we set the number of boundary points to $2p$. The classical FD-SBP operators
are constructed by considering degrees of freedom at the boundary nodes of the operator. Focusing on $p=2$, the operators have the following structure:
\begin{equation}
\Mclas :=\frac{2}{\Nxi}\diag\left(h_{1},\dots,h_{4},1,\dots,1,h_{4},\dots,h_{1} \right),
\end{equation}
and $\Dclas:=\left(\Mclas\right)^{-1}\Qclas$, where $\Qclas$ is given in Fig. \ref{fig:Q}.
\begin{figure}[!t]
\begin{center}
$\Qclas=
\begin{tikzpicture}[baseline=-0.5ex]
\draw [rounded corners,fill=gray, fill opacity=0.3] (-4.75,1.15)--(-4.75,3.25)--(-1.6,3.25)--(-1.6,1.15)--cycle;
\draw [rounded corners,fill=green, fill opacity=0.3] (-1.6,2.6)--(-1.6,1.15)--(0.6,1.15)--cycle;

\draw [rounded corners,fill=gray, fill opacity=0.3] (1.55,-3.25)--(1.55,-1.15)--(4.7,-1.15)--(4.7,-3.25)--cycle;
\draw [rounded corners,fill=green, fill opacity=0.3] (-0.5,-1.15)--(1.55,-1.15)--(1.55,-2.5)--cycle;

\draw [rounded corners, fill=mayablue, fill opacity=0.6] (-3.4,1.15)--(0.6,1.15)--(3.9,-1.15)--(-0.5,-1.15)--cycle;

\matrix[matrix of math nodes,left delimiter = {[},right delimiter = {]}] (m) {
-\frac{1}{2}&q_{12}&q_{13}&q_{14}&0&0&0&0&0&0&0&0\\
-q_{12}&0&q_{23}&q_{24}&0&0&0&0&0&0&0&0\\
-q_{13}&-q_{23}&0&q_{34}&-\frac{1}{12}&0&0&0&0&0&0&0\\
-q_{14}&-q_{24}&-q_{34}&0&\frac{2}{3}&-\frac{1}{12}&0&0&0&0&0&0\\
0&0&\frac{1}{12}&-\frac{2}{3}&0&\frac{2}{3}&-\frac{1}{12}&0&0&0&0&0\\
0&0&0&\frac{1}{12}&-\frac{2}{3}&0&\frac{2}{3}&-\frac{1}{12}&0&0&0&0\\
0&0&0&0&\frac{1}{12}&-\frac{2}{3}&0&\frac{2}{3}&-\frac{1}{12}&0&0&0\\
0&0&0&0&0&\frac{1}{12}&-\frac{2}{3}&0&\frac{2}{3}&-\frac{1}{12}&0&0\\
0&0&0&0&0&0&\frac{1}{12}&-\frac{2}{3}&0&q_{34}&q_{24}&q_{14}\\
0&0&0&0&0&0&0&\frac{1}{12}&-q_{34}&0&q_{23}&q_{13}\\
0&0&0&0&0&0&0&0&-q_{24}&-q_{23}&0&q_{12}\\
0&0&0&0&0&0&0&0&-q_{14}&-q_{13}&-q_{12}&\frac{1}{2}\\
};
\end{tikzpicture}$
\caption{Structure of degree 2 FD-SBP operator with free boundary parameters}
\label{fig:Q}
\end{center}
\end{figure}

Here the degree of $\Dclas$ on the interior nodes is $2p$, as we consider a central difference formula. For the classical FD-SBP operators the free coefficients are calculated, so that the matrix $\Dclas$ is of degree $p$. By doing this, the norm matrix $\Mclas$ is automatically of degree $2p-1$, see \cite{DCDRF2014,Hicken2013}.

We will consider degree preserving SBP operators $(\Mxione,\Dxione)$ with norm matrices
of degree $2p$. However, we construct our SBP operators in a similar fashion and consider an
interior stencil of degree $2p$ and boundary nodes with free parameters. Since the norm matrix needs
to be at least one degree higher than for the classical FD-SBP operators, one naturally needs more
free coefficients. These coefficients are obtained by increasing the number of boundary points at
least by one (to $2p+1$). For example, for $p=2$ the norm matrix $\Mxione$ is
\begin{equation}
\Mxione :=\frac{2}{\Nxi}\diag\left(h_{1},\dots,h_{5},1,\dots,1,h_{5},\dots,h_{1} \right),
\end{equation}
and the upper left corner of $\Qxione$ has the following structure:
\begin{center}
$\left(\Qxione\right)_{(1:5,1:7)}=
\begin{tikzpicture}[baseline=-0.5ex]
\matrix[matrix of math nodes,left delimiter = {[},right delimiter = {]}] (m) {
-\frac{1}{2}&q_{12}&q_{13}&q_{14}&q_{15}&0&0\\
-q_{12}&0&q_{23}&q_{24}&q_{25}&0&0\\
-q_{13}&-q_{23}&0&q_{34}&q_{35}&0&0\\
-q_{14}&-q_{24}&-q_{34}&0&q_{35}&-\frac{1}{12}&0\\
-q_{15}&-q_{25}&-q_{35}&-q_{45}&0&\frac{2}{3}&-\frac{1}{12}\\
};
\end{tikzpicture}.$
\end{center}
\textred{As for $\Qclas$ the degrees of freedom referring to the interior nodes are the same. Here,
we again consider a central difference formula, which is based on Taylor expansion. Due
to this stencil, we name the newly created operators \textit{degree preserving, element based finite difference} operators.
Note, that these operators are element based as in \cite{Fernandez2016corner}. By changing the number of nodes,
we must re-calculate the degrees of freedom at the boundary blocks. \cite{DCDRF2014}.}

Let $\bm{\xi}$ be a uniformly distributed set of nodes within the reference space $[-1,+1]$ and $\bm{\xi}_1=-1,~\bm{\xi}_{\Nxi}=+1$. The free coefficients are determined by solving
\begin{equation}\label{eq: degreeEqMod}
\begin{split}
\Qxione \bm{\xi}^k&=k\Mxione\bm{\xi}^{k-1}\qquad k=0,\dots,p,\\
\bm{1}\Tr\Mxione\bm{\xi}^{k}&=\frac{\bm{\xi}_{\Nxi}^{k+1}-\bm{\xi}_1^{k+1}}{k+1}\qquad k=0,\dots,2p.
\end{split}
\end{equation}
In comparison to classical FD-SBP operators, the degree preserving operators are
element based, in the sense that their coefficients explicitly depend on $\Nxi$. A disadvantage
of these operators is that the coefficients of the norm matrix are not necessarily positive
for an arbitrary choice of $\Nxi$.
In case of negative weights for a fixed value $\Nxi$, we increase the number of
boundary nodes to obtain more free coefficients for satisfying the property of a positive definite
norm matrix.

By solving \eqref{eq: degreeEqMod}, the SBP operators are not fully specified. As in
\cite{DCDRF2014}, the free parameters are chosen so that the truncation error is minimized. To do so we
define
\begin{equation}
\epsilon_{p+1}\coloneqq\left(\Mxione\right)^{-1}\Qxione\bm{\xi}^{p+1}-(p+1)\bm{\xi}^p,
\end{equation}
and then minimize
\begin{equation}\label{eq:trunc}
\bm{J}_e\coloneqq\epsilon_{p+1}\Tr\Mxione\epsilon_{p+1},
\end{equation}
under the assumption that all entries on the diagonal of $\Mxione$ are positive. Still, applying
this strategy does not necessarily specify all coefficients of $\Qxione$. Besides a small truncation
error, a desirable property is to have small coefficients within the operator to avoid round-off
errors. This is achieved by minimizing the sum of squares of the operator $\Qxione$ which is accomplished by
minimizing
\begin{equation}\label{eq:squares}
\bm{J}_Q\coloneqq\bm{1}\Tr\Qxione\circ\Qxione\bm{1},
\end{equation}
where $\circ$ is the Hadamard product, and we again ensure that all entries of $\Mxione$ are positive.

A description of the construction is provided in the following pseudo code:
\begin{center}
\begin{algorithm}[H]
 \KwIn{$\Nxi$ and $p$}
Set $bp:= 2p+1$\;
Solve \eqref{eq: degreeEqMod}\;
 \While{one or more entries of $\Mxione$ are negative }{
Set bp:= bp+1\;
Solve \eqref{eq: degreeEqMod} again\;
 }
Minimize $\bm{J}_e$ in \eqref{eq:trunc} under the \textblue{constraint} $h_i>0$ for $i=1,\dots,bp$\;
Minimize $\bm{J}_Q$ in \eqref{eq:squares} under the \textblue{constraint} $h_i>0$ for $i=1,\dots,bp$\;\
\caption{Construction of degree preserving operators}
\end{algorithm}
\end{center}

\textred{By following these steps we created SBP operators with a degree  $p$ differentiation matrix $\Dxione$ and a norm matrix $\Mxione$ with degree $2p$. Note, that as for existing SBP operators as in \cite{DCDRF2014}, a minimum number of nodes needs to be considered.}

We created the degree preserving operators using MAPLE.
\textred{The set of operators which \textblue{are} considered later in Sec. \ref{sec:results} were constructed with the MATLAB routines provided in the electronic supplementary material (ESM) of this article.}

\section{Proof of conservation, energy stability, and degree preservation}\label{sec:proof}
In this section, we prove conservation, energy stability, and the degree preserving property for the presented non-conforming discretization \eqref{eq:discA} and \eqref{eq:discB}. We consider degree preserving SBP operators where the degree of the differentiation matrices $\DxiL,\DxiR$ is $p$, and the degree of the norm matrices $\ML,\MR$ on the elements as well as the norm matrix $\MGam$ on the interface is equal to or greater than $2p$. Here, the choice of the nodes on the intermediate grid is completely arbitrary provided the norm matrix on the intermediate grid, $\MGam$, exists. For this reason we consider different nodal distributions on the intermediate grid in Sec. \ref{sec:results}.

\begin{thm}\label{thm:eins}
Given degree preserving SBP operators $\DxiA$ and $\DxiB$ of degree $p$, and projection operators of degree $\ge p$, then, the non-conforming discretizations of a left element
%
\begin{equation}\label{eq:discAThm}
\begin{split}
&\frac{\mr{d}\uA}{\mr{d}t}+\JA\lamxi\DxiA\uA+\JA\lametaA\DetaA\uA =\JA\HA^{-1}\left(\lamxi\ExiNL\uA-\fsL\right),
\end{split}
\end{equation}
a right element
%
\begin{equation}\label{eq:discBThm}
\begin{split}
&\frac{\mr{d}\uB}{\mr{d}t}+\JB\lamxi\DxiB\uB+\JB\lametaB\DetaB\uB =-\JB\HB^{-1}\left(\lamxi\ExioneR\uB-\fsR\right),
\end{split}
\end{equation}
and a single corresponding interface, where $\fsL$ and $\fsR$ are given in (\ref{eq:ustarLR}), the numerical approximation has the following properties:
\begin{enumerate}[label={\thethm.\arabic*}]
\item Discrete conservation, \textred{meaning the discrete integral of $\bm u$ is constant over time}. \label{Thmpartone}
\item Discrete energy stability. \label{Thmparttwo}
\item Discrete preservation of the degree $p$ for a non-conforming approximation. \label{Thmpartthree}
\end{enumerate}
\end{thm}
\begin{pf}
We prove the result in three parts.
\begin{proofpart}~~\ref{Thmpartone}:
Multiplying (\ref{eq:discAThm}) by $\fnc{J}_L \OneL\Tr\ML$ and ignoring the terms which are not related to the interface $\Gamma$, we have
\begin{equation}\label{eq:CProof}
\fnc{J}_L\OneL\Tr\ML\frac{\mr{d}\uL}{\mr{d}t}=-\lamxi \OneL\Tr\left(\SxiL+\frac{1}{2}\ExiNL\right)\uL+\OneL\Tr\left(\lamxi \ExiNL\uL-\fsL\right).
\end{equation}
Due to the SBP property and from the consistency of the derivative matrix $\Dxi \OneL=\bm{0}_{L}$, it holds that
\begin{equation}
\OneL\Tr\SxiL\uL=\frac{1}{2}\OneL\Tr\ExiNL\uL.
\end{equation}
Rearranging the right-hand side of \eqref{eq:CProof} yields
\begin{equation}
\begin{split}
\fnc{J}_L\OneL\Tr\ML\frac{\mr{d}\uL}{\mr{d}t}=&-\OneL\Tr\fsL,\\
=&-\frac{\lamxi}{2} \left(\OneL\Tr\ExiNL\uL+\OneL\Tr\RL\Tr\MGam\RR\uR \right)
-\frac{|\lamxi|\sigma}{2}\left(\OneL\Tr\RL\Tr\MGam\RR\uR-\OneL\Tr\RL\Tr\MGam\RL\uL  \right).
\end{split}
\end{equation}
From the properties of the projection operators \eqref{eq:RLRR} and the condition on the SATs \eqref{eq:dofinal}, it can be shown that
\begin{equation}\label{eq:ConL}
\begin{split}
\fnc{J}_L\OneL\Tr\ML\frac{\mr{d}\uL}{\mr{d}t}=&-\frac{\lamxi}{2} \left(\OneL\Tr\ExiNL\uL+\OneR\Tr\ExioneR\uR \right)
-\frac{|\lamxi|\sigma}{2}\left(\OneR\Tr\ExioneR\uR-\OneL\Tr\ExiNL\uL \right).
\end{split}
\end{equation}
Analogously, multiplying (\ref{eq:discBThm}) by $\fnc{J}_R \OneR\Tr\MR$ yields
\begin{equation}\label{eq:ConR}
\begin{split}
\fnc{J}_R\OneR\Tr\MR\frac{\mr{d}\uR}{\mr{d}t}=&\frac{\lamxi}{2} \left(\OneR\Tr\ExioneR\uR+\OneL\Tr\ExiNL\uL \right)+\frac{|\lamxi|\sigma}{2}\left(\OneR\Tr\ExioneR\uR-\OneL\Tr\ExiNL\uL  \right).
\end{split}
\end{equation}
Adding (\ref{eq:ConL}) and (\ref{eq:ConR}) we see that
\begin{equation}
\fnc{J}_L\OneL\Tr\ML\frac{\mr{d}\uL}{\mr{d}t}+\fnc{J}_R\OneR\Tr\MR\frac{\mr{d}\uR}{\mr{d}t}=0,
\end{equation}
which completes the proof.
\end{proofpart}
\begin{proofpart}~~\ref{Thmparttwo}:
To show energy stability we apply the discrete energy method. Multiplying (\ref{eq:discAThm}) by $\fnc{J}_L \uL\Tr\ML$ and ignoring the terms that are not related to the interface $\Gamma$, we obtain
\begin{equation}\label{eq:EProof}
\fnc{J}_L\uL\Tr\ML\frac{\mr{d}\uL}{\mr{d}t}=-\lamxi \uL\Tr\left(\SxiL+\frac{1}{2}\ExiNL\right)\uL+\uL\Tr\left(\lamxi\ExiNL\uL-\fsL\right).
\end{equation}
Since $\SxiL=-\SxiL\Tr$, we have
\begin{equation}\label{eq:newProp22}
\uL\Tr\SxiL\uL=0.
\end{equation}
Rearranging the right hand side of \eqref{eq:EProof}, \textred{including the numerical flux \eqref{eq:ustarLR}} and the property \eqref{eq:newProp22} yields
%
\begin{equation}\label{eq:EnL}
\begin{split}
\fnc{J}_L\uL\Tr\ML\frac{\mr{d}\uL}{\mr{d}t}=& \uL\Tr\left(\frac{\lamxi}{2}\ExiNL\uL-\fsL\right),\\
=&-\frac{\lamxi}{2} \uL\Tr\RL\Tr\MGam\RR\uR+\frac{|\lamxi|\sigma}{2}\left(\uL\Tr\RL\Tr\MGam\RR\uR-\uL\Tr\RL\Tr\MGam\RL\uL \right).
\end{split}
\end{equation}
For the discretization on the right element we multiply (\ref{eq:discBThm}) by $\fnc{J}_R\uR\Tr\MR$ to obtain
\begin{equation}\label{eq:EnR}
\begin{split}
\fnc{J}_R\uR\Tr\MR\frac{\mr{d}\uR}{\mr{d}t}=&\frac{\lamxi}{2} \uR\Tr\RR\Tr\MGam\RL\uL-\frac{|\lamxi|\sigma}{2}\left(\uR\Tr\RR\Tr\MGam\RR\uR-\uR\Tr\RR\Tr\MGam\RL\uL \right).
\end{split}
\end{equation}
Setting $\uGamL\coloneqq \RL\uL$ and $\uGamR\coloneqq \RR\uR$ and adding (\ref{eq:EnL}) and (\ref{eq:EnR}) we find that
\begin{equation}\label{eq:help2}
\begin{split}
\fnc{J}_L\uL\Tr\ML\frac{\mr{d}\uL}{\mr{d}t}+\fnc{J}_R\uR\Tr\MR\frac{\mr{d}\uR}{\mr{d}t}&=-\frac{|\lamxi|\sigma}{2}\left( \uGamR\Tr\MGam\uGamR-2\uGamR\Tr\MGam\uGamL+\uGamL\Tr\MGam\uGamL\right),\\
&=-\frac{|\lamxi|\sigma}{2}\begin{pmatrix}\uGamL\\ \uGamR\end{pmatrix}\Tr\underbrace{\left(\begin{pmatrix} +1 & -1 \\ -1 & +1\end{pmatrix}\right)\otimes\MGam}_{=:\mat{\tilde M}}\begin{pmatrix}\uGamL\\ \uGamR\end{pmatrix}\le 0.
\end{split}
\end{equation}
Here, $\mat{\tilde M}$ is positive semidefinite, since it is the tensor product of two positive semidefinite matrices. Since $|\lamxi|\sigma\ge 0$, the expression in \eqref{eq:help2} is negative, and the discretization remains energy stable.
\end{proofpart}
\begin{proofpart}~~\ref{Thmpartthree}:
Let $k=1,\dots,\nmin{p}$. For degree preservation the SAT on the corresponding interface must vanish for polynomials up to degree $p$. By including the polynomial $\bpkL$ in the SAT of (\ref{eq:discA}) we have
\begin{equation}\label{eq:ordpre}
SAT=\JL\ML^{-1}\left(\lamxi\ExiNL\bpkL-\fsL\right).
\end{equation}
The numerical flux is rearranged to be
\begin{equation}\label{eq:NumFluxOrd}
\begin{split}
\fsL&=\frac{\lamxi}{2}\left(\ExiNL\bpkL+\RL\Tr\MGam\RR\bpktR \right)-\frac{\sigma|\lamxi|}{2}\left( \RL\Tr\MGam\RR\bpktR-\RL\Tr\MGam\RL\bpkL \right),\\
&=\frac{\lamxi}{2}\left(\ExiNL\bpkL+\underbrace{\RL\Tr\MGam\bpkGamL}_{=\ExiNL\bpkL} \right)-\frac{\sigma|\lamxi|}{2}\left( \underbrace{\RL\Tr\MGam\bpkGamL-\RL\Tr\MGam\bpkGamL}_{=0} \right),\\
&=\lamxi\ExiNL\bpkL.
\end{split}
\end{equation}
Here, we have made use of the properties \eqref{eq:RLRR} and \eqref{eq:dofinal} of the projection operators. Substituting (\ref{eq:NumFluxOrd}) in (\ref{eq:ordpre}) we find
\begin{equation}
SAT=0.
\end{equation}
Hence, the resulting discretization is degree preserving for polynomials of degree $p$. We emphasize that if the SBP operators are not degree preserving, then projection operators that satisfy \eqref{eq:RLRR} and \eqref{eq:dofinal} cannot be constructed.\qed
\end{proofpart}
\end{pf}
%

\section{Numerical results}\label{sec:results}
In this section we demonstrate the convergence properties as well as the underlying theoretical findings from Thm. \ref{thm:eins} of the new non-conforming SBP-SAT scheme. To do so, we focus on the linear advection equation (\ref{LC}) in two dimensions on the domain $\Omega=[0,1]\times[0,1]$, where we assume a Cartesian mesh. The wave speeds in both directions are $\betax=\betay=1$. We enforce periodic boundary conditions to ensure that the energy is constant in the exact solution. Unless stated otherwise, the upwind SAT ($\sigma=1$) is chosen.
The time step is given by the CFL condition
\begin{equation}\label{eq:MaxCFL}
\Delta t:=CFL \frac{\min_{i}\{ \frac{\Delta x_i}{2},\frac{\Delta y_i}{2} \}}{\max_{i}\{N_i\}\max\{\alpha,\beta\}}
\end{equation}
where $\Delta x _i$ and $\Delta y _i$ denote the width in $x$- and $y$- direction of the $i$-th element, and $N_i$ denotes the number of nodes in one dimension of the $i$-th element. To integrate the approximation in time we use the five-stage, fourth-order low-storage Runge-Kutta method of Carpenter and Kennedy \cite{Kennedy1994}. Unless stated otherwise, we set $CFL=1$ and the final time to $T=0.1$. The error of the time integration does not affect our results.

For the convergence study, we consider the initial condition
$$\mathcal{U}(x,y,0)=2+\sin\left(2\pi x\right)+\cos\left(2\pi y\right).$$
The experimental order of convergence is determined by
\begin{equation*}
EOC_d=\frac{\log\left(\frac{L^2_d}{L^2_{d-1}}\right)  }{\log\left(\sqrt{\frac{DOFS_{d-1}}{DOFS_{d}}}\right)},
\end{equation*}
where $L^2_d$ denotes the error calculated with the norm matrix within \textred{all} elements at the $d$-th mesh level. This mesh level has $n_x=2^{d}$ elements in the $x$-direction and $n_y=2^{d}$ elements in the $y$-direction, which gives $4^{d}$ elements within the spatial domain. The $DOFS_d$ denotes the degrees of freedom in the domain. For our studies the number of nodes within an element remains constant, but the number of elements increases. This is a different approach than Kozdon and Wilcox \cite{Kozdon2016}, where for SBP finite difference operators the number of elements remains constant but the number of nodes within the elements increases.

\textred{This mesh refinement strategy is adapted from \cite{Gassner2013} as the newly derived SBP operators
are element based. For the newly derived SBP operators in Sec. \ref{sec:examples} their is no theoretical result concerning the
convergence order. However, we will numerically observe an order of $p+1$. }

The numerical results are divided into four components. First, in Sec. \ref{subsec:Comparison}, we demonstrate similar error and time step behavior for the new SBP operators compared to the finite-difference SBP operators, denoted by \textit{classical FD-SBP operators}, on conforming meshes. Next, Sec. \ref{subsec:Coupling Classic} demonstrates the order loss if we consider classical FD-SBP operators in the context of non-conforming problems. 
Sec. \ref{subsec:Coupling Degree} shows that the newly derived degree preserving SBP operators do not lose order in the convergence rate, for this same problem. Finally, Sec. \ref{subsec: Conservation and energy stability} verifies the proven conservation and energy stability from Thm. \ref{thm:eins} of the new scheme.
\subsection{Comparison of degree preserving SBP FD operators with classical FD-SBP operators on conforming meshes}\label{subsec:Comparison}
First, we consider a conforming grid to present the accuracy of the constructed degree preserving, element based finite difference SBP operators. We consider an SBP operator of degree $p=3$, where the degree of the norm matrix is $2p+1=7$. We compare this operator with the classical FD-SBP operator of degree $p=3$ adapted from \cite{Fernandez2014}, where the norm matrix is of degree $2p-1=5$. We note that the degree preserving and classical FD-SBP operators contain the same interior stencil. For comparison, we set the number of nodes in one dimension for both operators to be $22$.
For both operators we get a convergence rate of $p+1=4$ as shown in Tables \ref{tab:CompDegree} and \ref{tab:CompClassic}.
\begin{table}[!t]
\begin{center}
\begin{minipage}{0.4\textwidth}
\begin{center}
\begin{tabular}{c|c|c}
\hline
DOFS & $\Ltwo$ & EOC\\
\hline
1936&8.70E-05&\\
7744&5.80E-06&3.9\\
30976&3.03E-07&4.3\\
123904&1.98E-08&3.9\\
495616&1.24E-09&4.0\\
1982464&7.50E-11&4.0\\
\hline
\end{tabular}\\[0.1cm]
Maximum CFL number: 2.27
\caption{Experimental order of convergence and maximum CFL number for a \underline{degree preserving} FD operator with $p=3$.}
\label{tab:CompDegree}
\end{center}
\end{minipage}
\qquad
\begin{minipage}{0.4\textwidth}
\begin{center}
\begin{tabular}{c|c|c}
\hline
DOFS & $\Ltwo$ & EOC\\
\hline
1936&1.04E-04&\\
7744&7.71E-06&3.8\\
30976&5.44E-07&3.8\\
123904&3.27E-08&4.1\\
495616&2.02E-09&4.0\\
1982464&1.26E-10&4.0\\
\hline
\end{tabular}\\[0.1cm]
Maximum CFL number: 2.26
\caption{Experimental order of convergence and maximum CFL number for a \underline{classical} FD-SBP operator with $p=3$.}
\label{tab:CompClassic}
\end{center}
\end{minipage}
\end{center}
\end{table}

Comparing the errors, we see that the newly developed SBP operator has a smaller $\Ltwo$ error compared to the classical FD-SBP operator. The maximum allowable time step/maximum CFL number for the degree preserving operator is nearly identical to that for the classical FD-SBP operator.

For conforming meshes, these two operators have similar properties. However, since the degree of the norm matrix of the degree preserving operator is $\ge 2p$, we are able to couple these operators in a stable manner on non-conforming meshes without losing an order in the convergence rate, as demonstrated in Sec. \ref{subsec:Coupling Degree}. This is not possible for classical FD-SBP operators.

\subsection{Coupling classic SBP finite difference operators on non-conforming meshes}\label{subsec:Coupling Classic}
Here we demonstrate the order loss for non-conforming SBP schemes with classical FD-SBP operators by using the mesh refinement strategy as described above. We focus on the classical FD-SBP operator with $p=3$ as in Sec. \ref{subsec:Comparison}. To apply a non-conforming scheme we need to discretize $\Omega$ into elements with different nodal distributions. We divide the domain $\Omega$ into four subdomains:
\begin{equation*}
\begin{split}
\Omega_1 &= [0,0.5]\times[0,0.5],\\
\Omega_2 &= [0,0.5]\times[0.5,1],\\
\Omega_3 &= [0.5,1]\times[0,0.5],\\
\Omega_4 &= [0.5,1]\times[0.5,1].
\end{split}
\end{equation*}
In subdomains $\Omega_1$ and $\Omega_4$, we set all elements to have $N_1^2$ nodes ($N_1$ in both directions), whereas in $\Omega_2$ and $\Omega_3$ we set all elements to have $N_2^2$ nodes. For the simulations we considered two test cases with $N_1=22$ and $N_2=24$, as well as $N_1=22$ and $N_2=44$. For both configurations we have nodes that do not coincide on $\Gamma_1:=\{0.5\}\times[0,1]$ and $\Gamma_2:=[0,1]\times \{0.5\}$. For the degree $3$ operator we consider the projection operators derived by Kozdon and Wilcox \cite{Kozdon2016}. These operators are of degree $p-1=2$, so the non-conforming method is not degree preserving. Since the norm matrix of classical FD-SBP operators is $2p-1$ accurate, it is not possible to construct projection operators that are degree preserving and stable at the same time \cite{Nordstrom2015}. In Tables \ref{tab:Degree loss1} and \ref{tab:Degree loss2}, we see that these operators do not maintain an $EOC$ of $p+1$ or higher. In comparison, for the conforming case we achieve an $EOC$ of $p+1$. In the next section we will see, that by considering degree preserving SBP operators the non-conforming scheme maintains an $EOC$ of $p+1$ of higher.
\begin{table}[!t]
\begin{center}
\begin{minipage}{0.4\textwidth}
\begin{center}
$N_1=22$ and $N_2=24$\\
\begin{tabular}{c|c|c}
\hline
DOFS & $\Ltwo$ & EOC\\
\hline
2120 & 4.93E-04 & \\
8480 & 5.45E-05 & 3.2\\
33920 & 5.90E-06 & 3.2\\
135680 & 6.13E-07 & 3.3\\
542720 & 6.72E-08 & 3.2\\
2170880 & 7.54E-09 & 3.2\\
\hline
\end{tabular}\\[0.1cm]
Maximum CFL number: 2.24
\caption{Non-conforming method with classical FD-SBP operators and \underline{$N_2=N_1+2$}.}
\label{tab:Degree loss1}
\end{center}
\end{minipage}
\qquad
\begin{minipage}{0.4\textwidth}
\begin{center}
$N_1=22$ and $N_2=44$\\
\begin{tabular}{c|c|c}
\hline
DOFS & $\Ltwo$ & EOC\\
\hline
4840&4.27E-04&\\
19360&4.63E-05&3.2\\
77440&5.04E-06&3.2\\
309760&5.35E-07&3.2\\
1239040&5.96E-08&3.2\\
4956160& 6.82E-09 &3.1 \\
\hline
\end{tabular}\\[0.1cm]
Maximum CFL number: 2.16
\caption{Non-conforming method with classical FD-SBP operators and \underline{$N_2=2N_1$}.}
\label{tab:Degree loss2}
\end{center}
\end{minipage}
\end{center}
\end{table}\\
\subsection{Coupling of degree preserving SBP operators on non-conforming meshes}\label{subsec:Coupling Degree}
To demonstrate the efficiency of the two dimensional degree preserving SBP-SAT scheme, we construct the SATs using tensor products, as briefly described in Appendix \ref{sec:Projection}. By doing so we can directly compare our non-conforming method with already existing coupling methods \cite{Kozdon2016,Mattsson2010b}, since these methods are also tensor product based.

For this numerical test we couple the new degree preserving SBP operators. To do so we divide $\Omega$ in the same way as in Sec. \ref{subsec:Coupling Classic}. On each of the subdomains we consider the same nodal distribution, so the discretizations do not differ. For $\Gamma_1$ and $\Gamma_2$ we construct tensor-product based projection operators, as discussed in Appendix \ref{sec:Projection}.
For the numerical results we consider three choices for the intermediate grid nodes with a corresponding surface-norm matrix:
\begin{enumerate}
\item The nodes and norm matrix in one dimension from elements within $\Omega_1$ and $\Omega_4$, Fig. \ref{fig:Left}.
\item The nodes and norm matrix in one dimension from elements within $\Omega_2$ and $\Omega_3$, Fig. \ref{fig:Right}.
\item Twice the amount of nodes as on the most dense element in one dimension with the norm matrix adapted from the SBP operator, Fig. \ref{fig:twice}.
\end{enumerate}
\begin{figure}[!t]
\begin{center}
\begin{tikzpicture}[scale=0.5]
\draw[very thick] (-5,-2) -- (-1,-2) -- (-1,2) -- (-5,2) -- cycle;
\draw[very thick] (1,-2) -- (5,-2) -- (5,2) -- (1,2) -- cycle;
\draw[very thick] (0,3) -- (0,-3) ;
\color[rgb]{0,0,1}
\draw[] (-4,2) -- (-4,-2) ;
\draw[] (-3,2) -- (-3,-2) ;
\draw[] (-2,2) -- (-2,-2) ;
\draw[] (-5,1) -- (-1,1) ;
\draw[] (-5,0) -- (-1,0) ;
\draw[] (-5,-1) -- (-1,-1) ;
\node at (-1,-2) {X};
\node at (-1,-1) {X};
\node at (-1,0) {X};
\node at (-1,1) {X};
\node at (-1,2) {X};
\color[rgb]{1,0,0}
\draw[] (2.33333,2) -- (2.33333,-2) ;
\draw[] (3.66666,2) -- (3.66666,-2) ;
\draw[] (1,0.66666) -- (5,0.66666) ;
\draw[] (1,-0.66666) -- (5,-0.66666) ;
\node at (1,-2) {X};
\node at (1,-0.66666) {X};
\node at (1,0.66666) {X};
\node at (1,2) {X};
\color[rgb]{0,0,1}
\color[rgb]{0,0,1}
\node at (0,-2) {X};
\node at (0,-1) {X};
\node at (0,0) {X};
\node at (0,1) {X};
\node at (0,2) {X};
\end{tikzpicture}
\end{center}
		\caption{Interface grid nodes coincide with the left element.}
		\label{fig:Left}
\begin{center}
\begin{tikzpicture}[scale=0.5]
\draw[very thick] (-5,-2) -- (-1,-2) -- (-1,2) -- (-5,2) -- cycle;
\draw[very thick] (1,-2) -- (5,-2) -- (5,2) -- (1,2) -- cycle;
\draw[very thick] (0,3) -- (0,-3) ;
\color[rgb]{0,0,1}
\draw[] (-4,2) -- (-4,-2) ;
\draw[] (-3,2) -- (-3,-2) ;
\draw[] (-2,2) -- (-2,-2) ;
\draw[] (-5,1) -- (-1,1) ;
\draw[] (-5,0) -- (-1,0) ;
\draw[] (-5,-1) -- (-1,-1) ;
\node at (-1,-2) {X};
\node at (-1,-1) {X};
\node at (-1,0) {X};
\node at (-1,1) {X};
\node at (-1,2) {X};
\color[rgb]{1,0,0}
\draw[] (2.33333,2) -- (2.33333,-2) ;
\draw[] (3.66666,2) -- (3.66666,-2) ;
\draw[] (1,0.66666) -- (5,0.66666) ;
\draw[] (1,-0.66666) -- (5,-0.66666) ;
\node at (1,-2) {X};
\node at (1,-0.66666) {X};
\node at (1,0.66666) {X};
\node at (1,2) {X};

\color[rgb]{1,0,0}
\node at (0,-2) {X};
\node at (0,-0.66666) {X};
\node at (0,0.66666) {X};
\node at (0,2) {X};
\end{tikzpicture}
\end{center}
		\caption{Interface grid nodes coincide with the right element.}
		\label{fig:Right}
		\begin{center}
\begin{tikzpicture}[scale=0.5]
\draw[very thick] (-5,-2) -- (-1,-2) -- (-1,2) -- (-5,2) -- cycle;
\draw[very thick] (1,-2) -- (5,-2) -- (5,2) -- (1,2) -- cycle;
\draw[very thick] (0,3) -- (0,-3) ;
\color[rgb]{0,0,1}
\draw[] (-4,2) -- (-4,-2) ;
\draw[] (-3,2) -- (-3,-2) ;
\draw[] (-2,2) -- (-2,-2) ;
\draw[] (-5,1) -- (-1,1) ;
\draw[] (-5,0) -- (-1,0) ;
\draw[] (-5,-1) -- (-1,-1) ;
\node at (-1,-2) {X};
\node at (-1,-1) {X};
\node at (-1,0) {X};
\node at (-1,1) {X};
\node at (-1,2) {X};
\color[rgb]{1,0,0}
\draw[] (2.33333,2) -- (2.33333,-2) ;
\draw[] (3.66666,2) -- (3.66666,-2) ;
\draw[] (1,0.66666) -- (5,0.66666) ;
\draw[] (1,-0.66666) -- (5,-0.66666) ;
\node at (1,-2) {X};
\node at (1,-0.66666) {X};
\node at (1,0.66666) {X};
\node at (1,2) {X};
\color[rgb]{0,1,0}
\node at (0,-2) {X};
\node at (0,-1.5) {X};
\node at (0,-1) {X};
\node at (0,-0.75) {X};
\node at (0,-0.25) {X};
\node at (0,0.25) {X};
\node at (0,0.75) {X};
\node at (0,1) {X};
\node at (0,1.5) {X};
\node at (0,2) {X};
\end{tikzpicture}
\end{center}
		\caption{Interface has twice the amount of nodes as on the left element.}
		\label{fig:twice}
\end{figure}
For each set of nodes, the norm matrix on the intermediate grid is \textblue{of degree} $>2p$, so it is possible to create the appropriate projection matrices. However, no matter how we choose the intermediate grid, we achieve a convergence rate $\ge p+1=4$, as shown in Tables \ref{tab:1}-\ref{tab:6}.
\begin{table}[!t]
\begin{center}
\begin{minipage}{0.4\textwidth}
\begin{center}
$N_1=22$ and $N_2=24$\\
\begin{tabular}{c|c|c}
\hline
DOFS & $\Ltwo$ & EOC\\
\hline
2120&4.23E-03&\\
8480&3.06E-04&3.8\\
33920&1.97E-05&4.0\\
135680&1.15E-06&4.1\\
542720&6.52E-08&4.1\\
2170880&3.62E-09&4.2\\
\hline
\end{tabular}\\[0.1cm]
Maximum CFL number: 2.25
\caption{Degree preserving SBP operators: Intermediate grid with $22$ nodes and \underline{$N_2=N_1+2$}.}
\label{tab:1}
\end{center}
\end{minipage}
\qquad
\begin{minipage}{0.4\textwidth}
\begin{center}
$N_1=22$ and $N_2=44$\\
\begin{tabular}{c|c|c}
\hline
DOFS & $\Ltwo$ & EOC\\
\hline
4840&2.83E-03&\\
19360&2.01E-04&3.8\\
77440&1.21E-05&4.1\\
309760&6.50E-07&4.2\\
1239040&3.20E-08&4.3\\
4956160& 1.59E-09&4.3\\
\hline
\end{tabular}\\[0.1cm]
Maximum CFL number: 2.17
\caption{Degree preserving SBP operators: Intermediate grid with $22$ nodes and \underline{$N_2=2N_1$}.}
\label{tab:2}
\end{center}
\end{minipage}
\end{center}
\end{table}

\begin{table}[!t]
\begin{center}
\begin{minipage}{0.4\textwidth}
\begin{center}
$N_1=22$ and $N_2=24$\\
\begin{tabular}{c|c|c}
\hline
DOFS & $\Ltwo$ & EOC\\
\hline
2120&5.18E-03&\\
8480&3.56E-04&3.9\\
33920&2.19E-05&4.0\\
135680&1.27E-06&4.1\\
542720&7.28E-08&4.1\\
2170880&4.06E-09&4.2\\
\hline
\end{tabular}\\[0.1cm]
Maximum CFL number: 2.25
\caption{Degree preserving SBP operators: Intermediate grid with $24$ nodes and \underline{$N_2=N_1+2$}.}
\label{tab:3}
\end{center}
\end{minipage}
\qquad
\begin{minipage}{0.4\textwidth}
\begin{center}
$N_1=22$ and $N_2=44$\\
\begin{tabular}{c|c|c}
\hline
DOFS & $\Ltwo$ & EOC\\
\hline
4840&1.91E-03&\\
19360&1.31E-04&3.9\\
77440&8.30E-06&4.0\\
309760&4.77E-07&4.1\\
1239040&2.67E-08&4.2\\
4956160&1.51E-09&4.1\\
\hline
\end{tabular}\\[0.1cm]
Maximum CFL number: 2.17
\caption{Degree preserving SBP operators: Intermediate grid with $44$ nodes and \underline{$N_2=2N_1$}.}
\label{tab:4}
\end{center}
\end{minipage}
\end{center}
\end{table}

\begin{table}[!t]
\begin{center}
\begin{minipage}{0.4\textwidth}
\begin{center}
$N_1=22$ and $N_2=24$\\
\begin{tabular}{c|c|c}
\hline
DOFS & $\Ltwo$ & EOC\\
\hline
2120&1.05E-03&\\
8480&7.05E-05&3.9\\
33920&4.17E-06&4.1\\
135680&2.28E-07&4.2\\
542720&1.29E-08&4.1\\
2170880&7.16E-10&4.2\\
\hline
\end{tabular}\\[0.1cm]
Maximum CFL number: 2.25
\caption{Degree preserving SBP operators: Intermediate grid with $48$ nodes and \underline{$N_2=N_1+2$}.}
\label{tab:5}
\end{center}
\end{minipage}
\qquad
\begin{minipage}{0.4\textwidth}
\begin{center}
$N_1=22$ and $N_2=44$\\
\begin{tabular}{c|c|c}
\hline
DOFS & $\Ltwo$ & EOC\\
\hline
4840& 6.56E-04&\\
19360& 4.10E-05&4.0\\
77440& 2.25E-06&4.2\\
309760& 1.21E-07&4.2\\
1239040& 6.93E-09&4.1\\
4956160& 4.07E-10&4.1\\
\hline
\end{tabular}\\[0.1cm]
Maximum CFL number: 2.17
\caption{Degree preserving SBP operators: Intermediate grid with $88$ nodes and \underline{$N_2=2N_1$}.}
\label{tab:6}
\end{center}
\end{minipage}
\end{center}
\end{table}

\begin{figure}[!t]
\begin{center}
	\begin{minipage}{0.4\textwidth}
\begin{flushleft}
		\vspace{0.1cm}
	  \includegraphics[width=\textwidth]{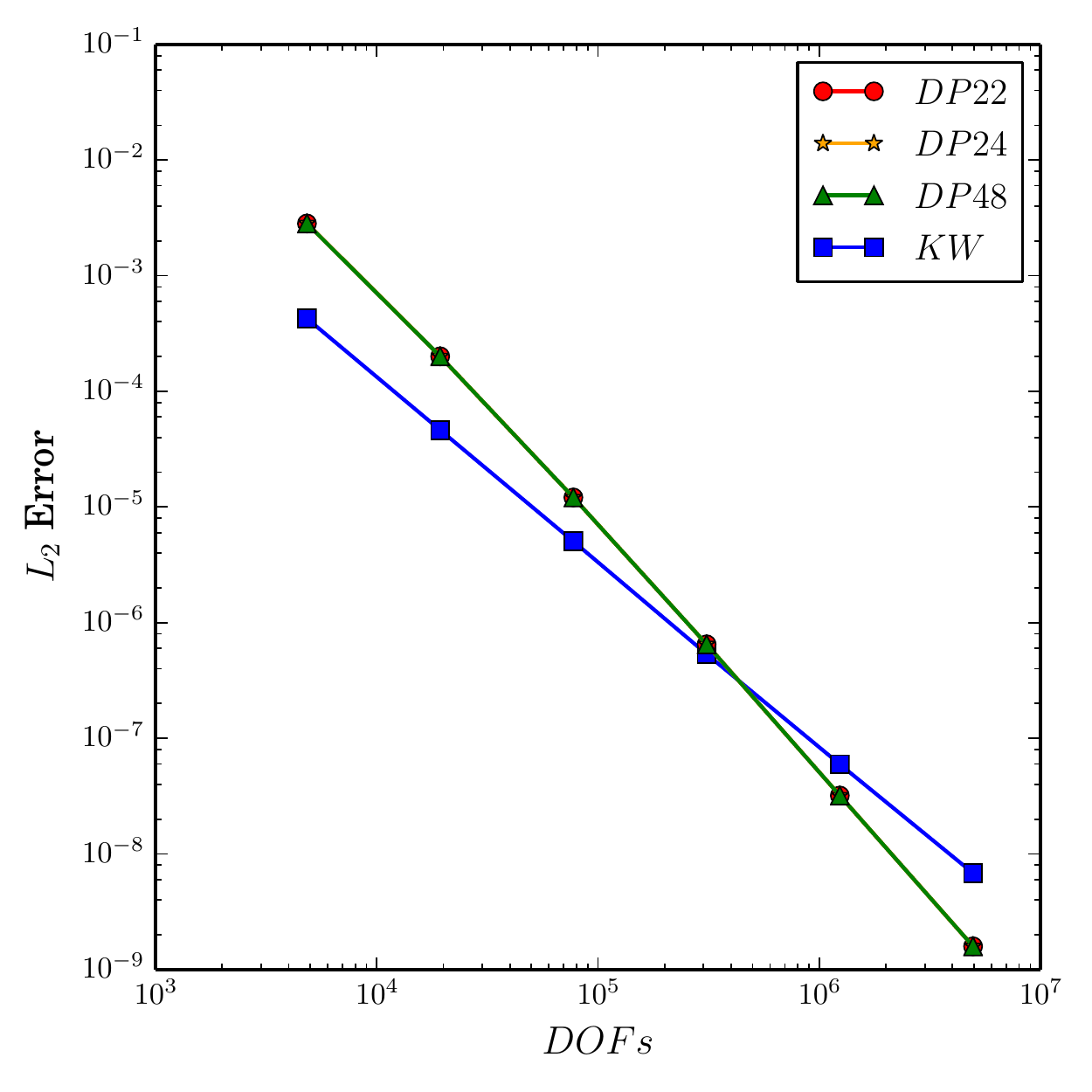}
		\caption{$L_2$ error for $N_2=N_1+2$ and different intermediate grids.}
		\label{fig: N1p2}
\end{flushleft}
	\end{minipage}
	\qquad
	\begin{minipage}{0.4\textwidth}
		\vspace{0.1cm}
	  \includegraphics[width=\textwidth]{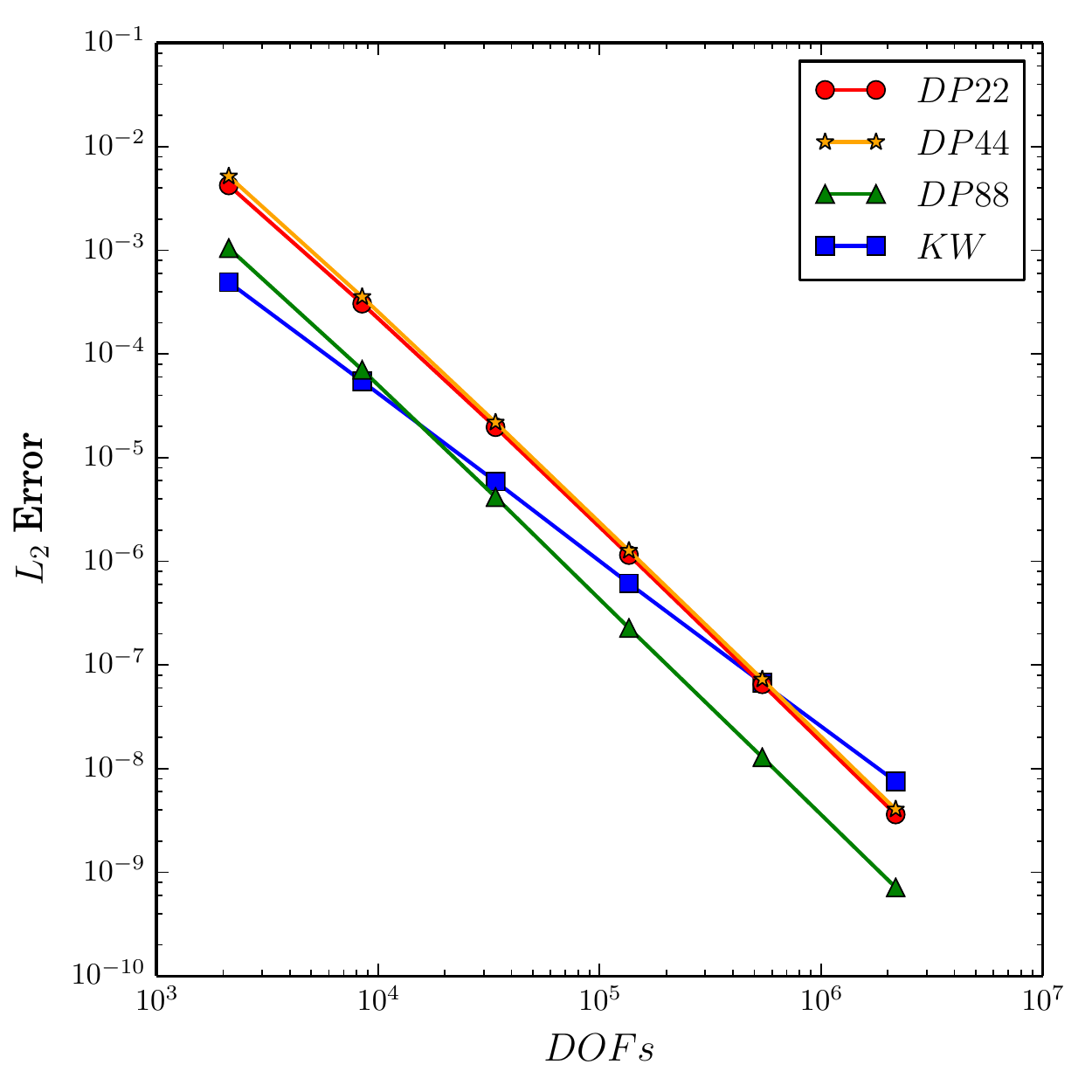}
		\caption{$L_2$ error for $N_2=2N_1$ and different intermediate grids.}
		\label{fig: 2N1}
	\end{minipage}
	\end{center}
\end{figure}

Considering the mesh configuration where the DOFs are the same, there is no change in the maximum CFL number. However, the $\Ltwo$ error is smaller if the intermediate grid contains more nodes, see Figs. \ref{fig: N1p2} and \ref{fig: 2N1}. Within these error plots, \emph{DPxx} denotes the simulation with degree preserving operators where the intermediate grid consists of \emph{xx} nodes, and \emph{KW} denotes the simulation with classical FD-SBP operators and the projection operators of Kozdon and Wilcox \cite{Kozdon2016}.


Note that the degree of the differentiation matrix for all simulations is the same ($p=3$).
The EOC of about one order higher is due to the fact that we can construct projection operators of one degree higher than with classical FD-SBP operators as a result of using degree preserving operators.
\subsection{Numerical verification of conservation and energy stability}\label{subsec: Conservation and energy stability}
As was shown in Thm. \ref{thm:eins}, the new scheme is conservative and energy stable. To verify this result we use a more complicated mesh distribution. We divide the domain $\Omega$ into eight elements in the $x$-direction and $y$-direction, resulting in a mesh of $64$ total elements. We divide the mesh distribution into a checker board pattern of different node distributions. Imagine that the tiles of the checker board are elements. In all black tiles we consider the degree preserving operator with $p=3$ and $22$ nodes in each direction. For all white tiles we consider the same operator with $24$ nodes in each direction. By distributing the mesh in this way, there is no interface in the interior of the domain where the nodes are conforming. So a projection  must be taken into account at each interface. With this mesh distribution we consider two simulations, one using a upwind flux $\sigma=1$ and one using a central flux $\sigma=0$.

\textred{We verify the stability results by writing the fully coupled system as a linear equation
\begin{equation}
\frac{d\bm u}{dt}=\mat A\bm u.
\end{equation}
By computing the eigenvalues of A we get the spectrum plot in Fig. \ref{fig: EigUp} and \ref{fig: EigCen}.\\
\begin{figure}[!t]
\begin{center}
	\begin{minipage}{0.4\textwidth}
\begin{flushleft}
		\vspace{0.1cm}
	  \includegraphics[width=\textwidth]{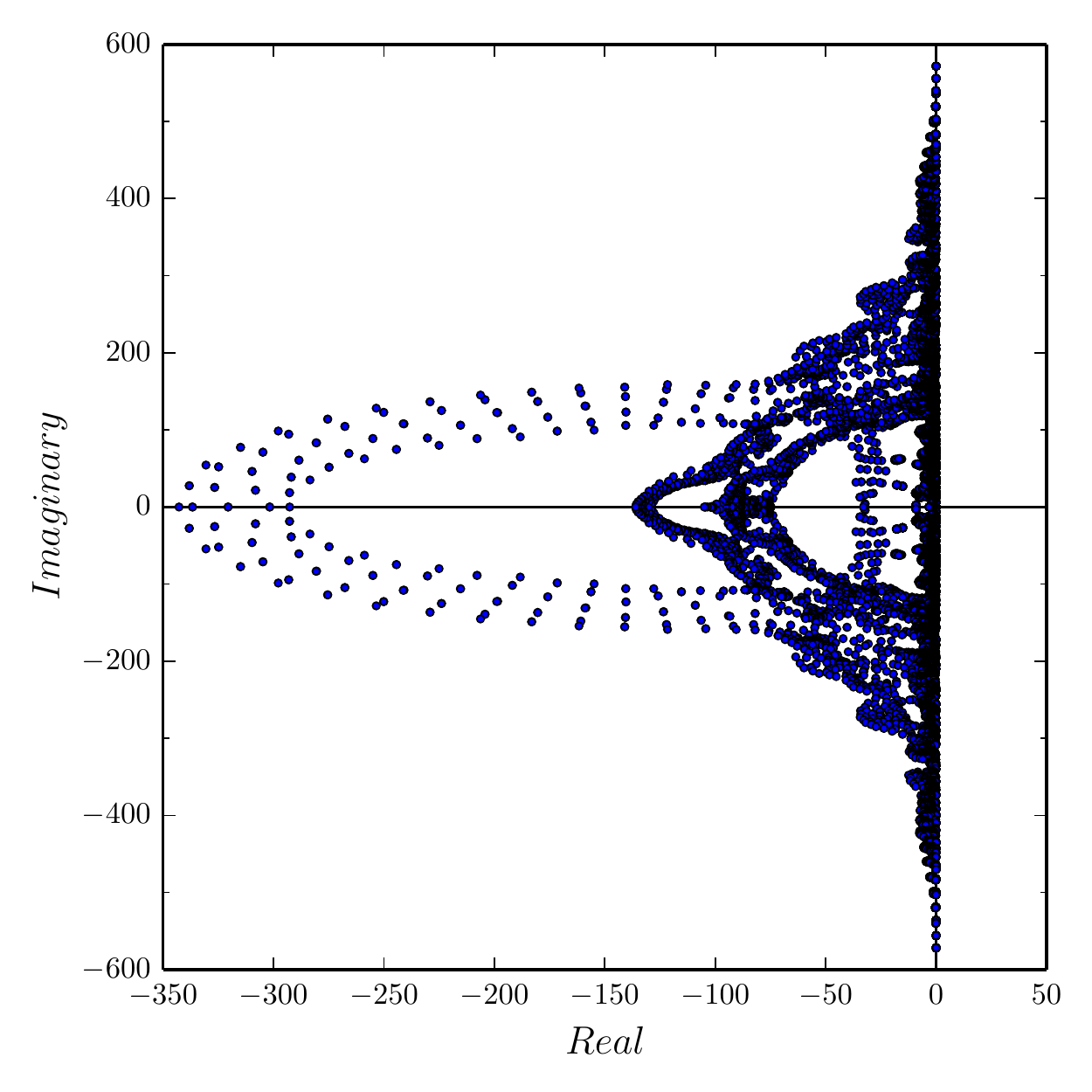}
		\caption{Non-conforming scheme with an upwind flux. All real parts of the eigenvalues are non-positive which indicates energy stability.}
		\label{fig: EigUp}
\end{flushleft}
	\end{minipage}
	\qquad
	\begin{minipage}{0.4\textwidth}
		\vspace{0.1cm}
	  \includegraphics[width=\textwidth]{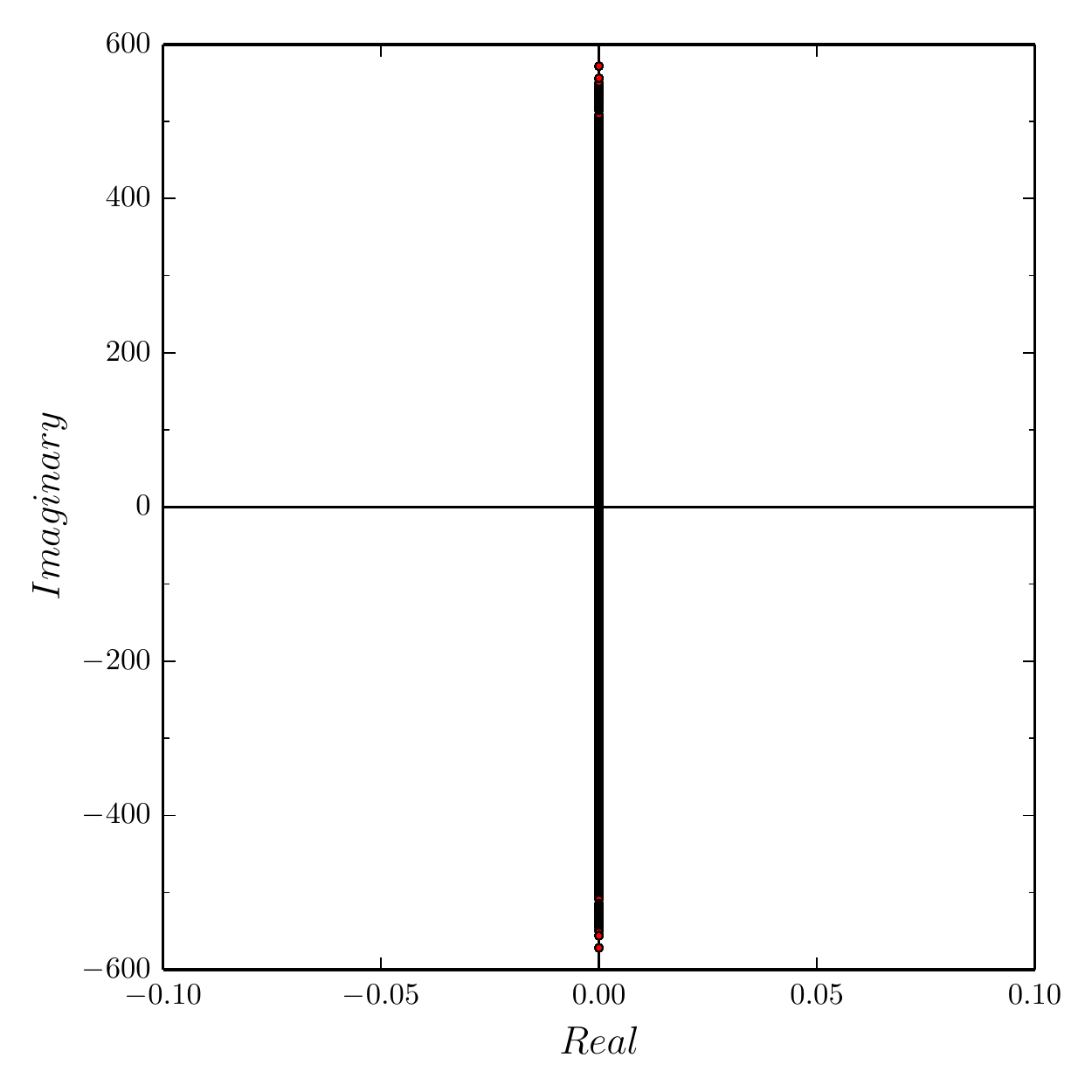}
		\caption{Non-conforming scheme with a central flux. All real parts of the eigenvalues are zero (besides round-off errors) which indicates energy conservation.}
		\label{fig: EigCen}
	\end{minipage}
	\end{center}
\end{figure}\\
For all eigenvalues the real parts are non-positive demonstrating an energy stable scheme. Furthermore, when considering a central flux all eigenvalues are purely imaginary which leads to energy conservative scheme. }

\textred{Besides the eigenvalues spectrum we analyze the behaviour of energy over time. Therefore, we consider simulations with $T=10$ and $CFL=0.25$.} Throughout the simulations, we measure conservation and energy in terms of
\begin{equation}
\begin{split}
conservation\; metric:=&\Bigg|\sum_{i=1}^K\fnc{J}_i\bm{1}\Tr\mat{H}_i\bm{u}_i(T)-\sum_{i=1}^K\fnc{J}_i\bm{1}\Tr\mat{H}_i\bm{u}_i(0)\Bigg|,\\
energy:=&\sum_{i=1}^K\fnc{J}_i\bm{u}_i\Tr\mat{H}_i\bm{u}_i,
\end{split}
\end{equation}
where $K$ denotes the number of elements and $\fnc{J}_i,\mat{H}_i, \bm{u}_i$ denote the mapping term, norm matrix and solution of the corresponding element respectively.
As initial condition we consider the discontinuous function, as this is a demanding test of the conservative properties of the scheme,
\begin{equation}
\mathcal{U}(x,y,t)=\begin{cases}
     3  & \text{for } x\le 0.3 \\
     1  &\text{for } x>0.3
   \end{cases},
\end{equation}
with periodic boundary conditions. To show that the new scheme is conservative and energy stable for nearly arbitrary nodal distributions on the intermediate grid, we consider the minimum  number of nodes case where $\MGam$ has the minimum nodes required to be at least of degree $2p$. To do so we consider an intermediate grid with four Legendre-Gauss nodes with corresponding weights for the diagonal entries of $\MGam$.
\begin{figure}[!t]
\begin{center}
	\begin{minipage}{0.4\textwidth}
\begin{flushleft}
		\vspace{0.1cm}
	  \includegraphics[width=\textwidth]{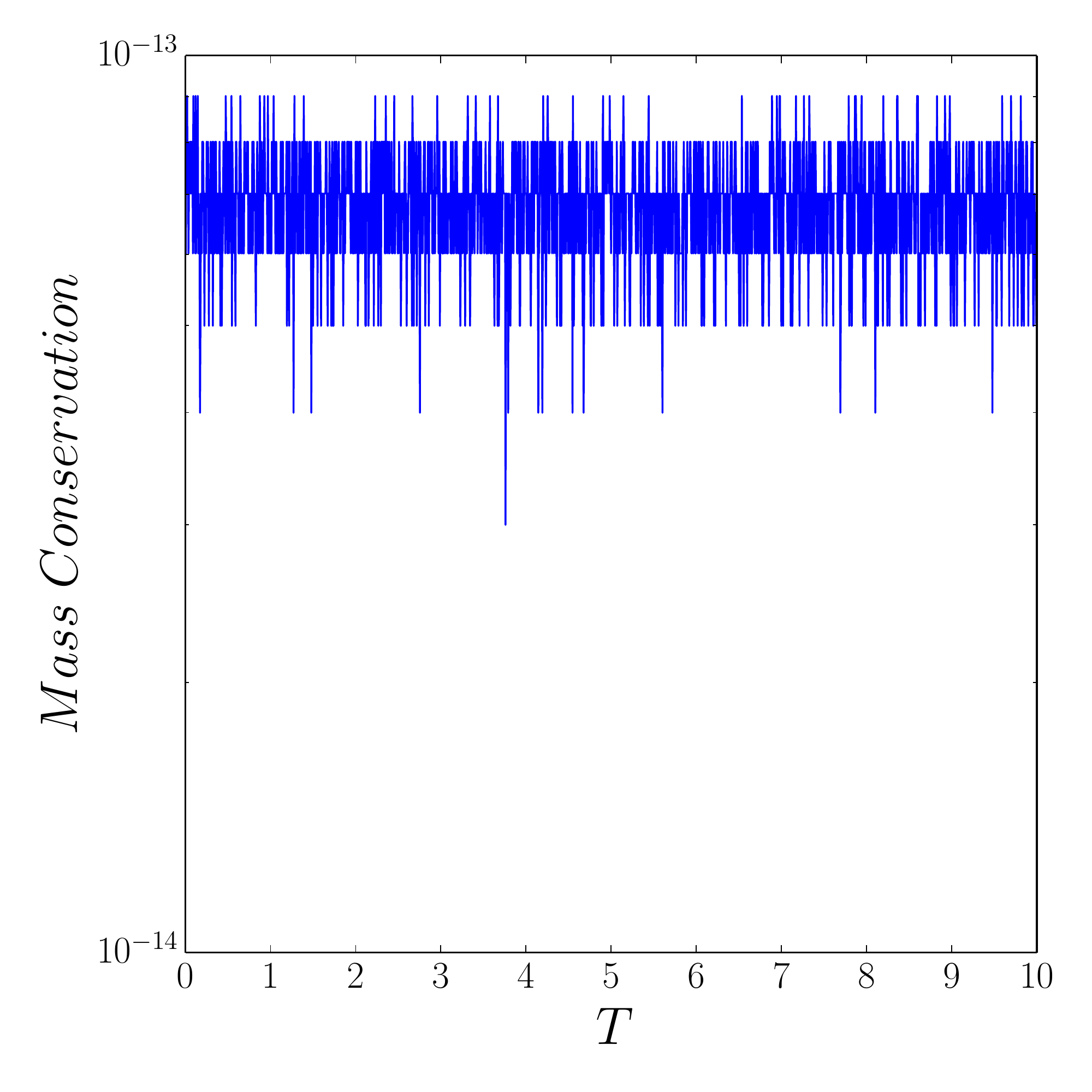}
		\caption{A semilog plot that demonstrates the evolution of the absolute difference of the initial total mass and current total mass of the solution. This demonstrates numerically that the degree preserving scheme is globally conservative.}
		\label{fig: Conservation}
\end{flushleft}
	\end{minipage}
	\qquad
	\begin{minipage}{0.4\textwidth}
		\vspace{0.1cm}
	  \includegraphics[width=\textwidth]{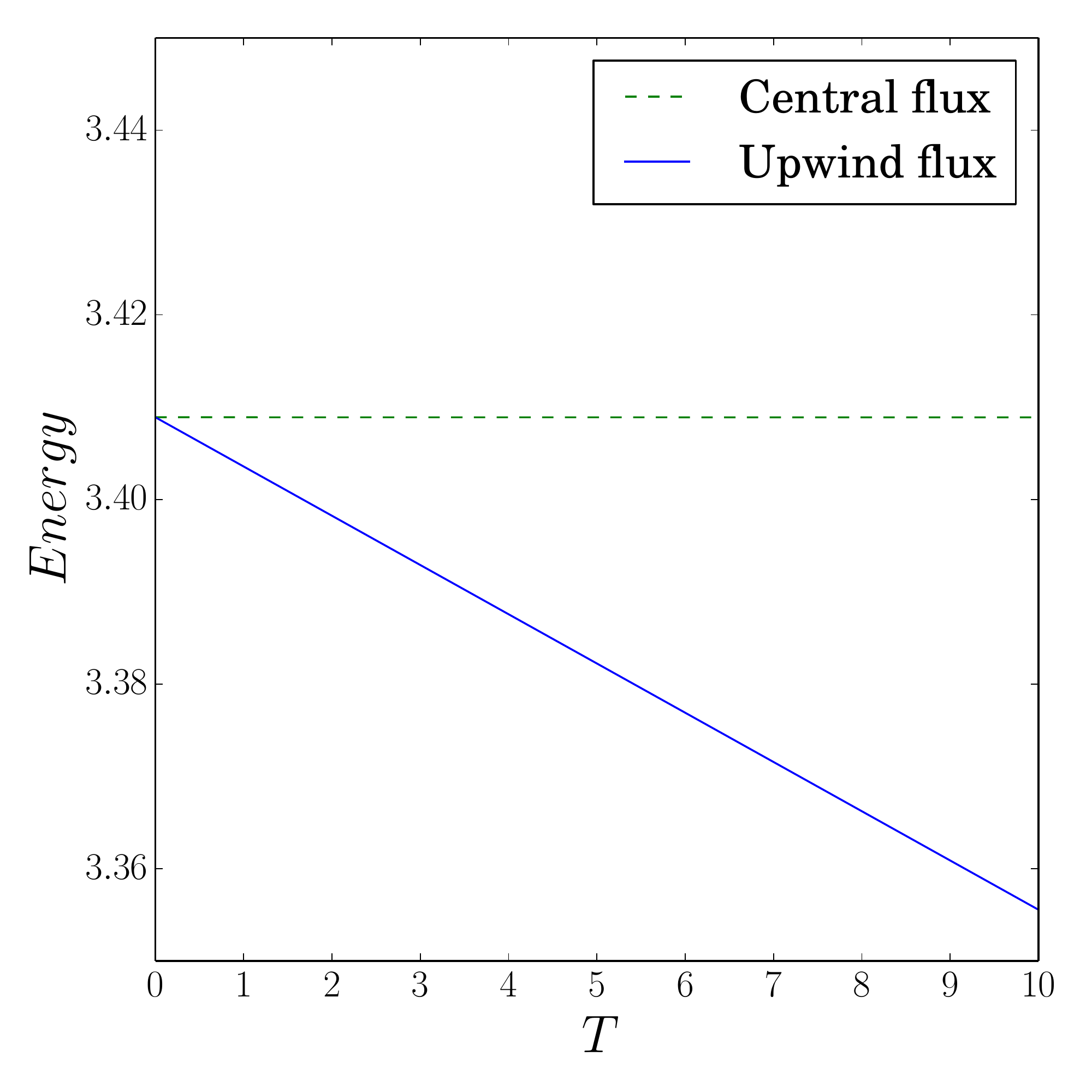}
		\caption{Evolution of the energy of the solution for a computation coupled with SATs using either an upwind or centered flux. We see that energy is conserved for the central scheme and energy decays for the upwind scheme as predicted by the theory.}
		\label{fig: Energy}
	\end{minipage}
	\end{center}
\end{figure}\\

As we can see in Fig.~\ref{fig: Conservation}, conservation holds. By using an upwind SAT, energy is dissipated, whereas using a central SAT the energy remains constant in time, as shown in Fig. \ref{fig: Energy}. For this test case the intermediate grid has four Gauss nodes, whereas there are many more nodes on the elements. This numerically verifies conservation and energy stability for this particular test case. However, as proven in Sec. \ref{sec:proof}, this behavior occurs for all kinds of mesh distributions.

\begin{rmk}
By considering intermediate grids with minimum number of Gauss nodes and applying the same test as in Sec. \ref{subsec:Coupling Degree}, we also maintain an EOC higher than $p+1$ as shown in Tables \ref{tab:G1} and \ref{tab:G2}.

\begin{table}[!t]
\begin{center}
\begin{minipage}{0.4\textwidth}
\begin{center}
$N_1=22$ and $N_2=24$\\
\begin{tabular}{c|c|c}
\hline
DOFS & $\Ltwo$ & EOC\\
\hline
2120&1.20E-02&\\
8480&7.73E-04&4.0\\
33920&4.66E-05&4.1\\
135680&2.53E-06&4.2\\
542720&1.36E-07&4.2\\
2170880&7.33E-09&4.2\\
\hline
\end{tabular}\\[0.1cm]
Maximum CFL number: 2.24
\caption{Degree preserving SBP operators: Intermediate grid with $4$ Gauss nodes and \underline{$N_2=N_1+2$}.}
\label{tab:G1}
\end{center}
\end{minipage}
\qquad
\begin{minipage}{0.4\textwidth}
\begin{center}
$N_1=22$ and $N_2=44$\\
\begin{tabular}{c|c|c}
\hline
DOFS & $\Ltwo$ & EOC\\
\hline
4840& 1.21E-02&\\
19360& 7.65E-04&4.0\\
77440& 4.59E-05&4.1\\
309760& 2.51E-06&4.2\\
1239040& 1.37E-07&4.2\\
4956160& 7.67E-09&4.2\\
\hline
\end{tabular}\\[0.1cm]
Maximum CFL number: 2.16
\caption{Degree preserving SBP operators: Intermediate grid with $4$ Gauss nodes and \underline{$N_2=2N_1$}.}
\label{tab:G2}
\end{center}
\end{minipage}
\end{center}
\end{table}

However, the $\Ltwo$ error is larger than the presented ones in Sec. \ref{subsec:Coupling Degree}. This underlines the assumption that the number of nodes on the intermediate grid is related to the overall $\Ltwo$ error.
\end{rmk}

\section{Conclusions}\label{sec:conclusions}
In this paper we have derived a non-conforming SBP scheme built from multi-dimensional SATs. When
tensor products are used, such non-conforming schemes have already been derived as in
\cite{Kozdon2016,Mattsson2010b}, where classical FD-SBP finite difference operators are considered.
However, due to the use of classical FD-SBP finite difference operators, it is not possible to
construct projection operators between non-conforming interfaces that preserve the degree of the
derivative matrix while maintaining energy stability, as proven by Lundquist and Nordstr\"om and given in Thm. \ref{thrm:tpm1}.
This results from the fact that the norm matrix, $\mat{H}^{(1\textrm{D})}$, is of degree $2p-1$ and therefore projection operators
can not be constructed that are of degree $\ge p$. Having identified the norm matrix as the factor
impeding the construction of appropriate projection operators, we focused on the derivation of degree preserving operators in one
dimension. To do so we increased the degree of the norm matrix of the degree preserving SBP
operators to be higher than that used in classical FD-SBP operators.

The new degree preserving SBP operators allowed us to construct a non-conforming, degree preserving, conservative, and energy stable method. These properties of the scheme were verified analytically as well as numerically. The main idea of the method is to project the solution from each element to a set of nodes on an intermediate grid on the interface between elements and evaluate the SAT term on these intermediate nodes. On the intermediate grid the nodal distribution is arbitrary, as long as a quadrature rule of a minimum degree exists. In the numerical results, we also considered different intermediate grids.

By comparing our newly developed degree preserving scheme with already existing non-conforming methods, we achieved experimental orders of convergence of one order higher; this is likely a result of the fact that the new scheme is globally of degree $p$.
\section*{Acknowledgments}

The work of Lucas Friedrich, Andrew Winters and Gregor Gassner was funded by the Deutsche Forschungsgemeinschaft (DFG) grant TA 2160/1-1.

The work of Jason Hicken was partially funded by the Air Force Office of Scientific Research Award FA9550-15-1-0242 under Dr. Jean-Luc Cambier.

This work was partially performed on the Cologne High Efficiency Operating Platform for Sciences (CHEOPS) at the Regionales Rechenzentrum K\"oln (RRZK).

\bibliographystyle{siam}
\bibliography{Bib}

\appendix


\section{Projection operators}\label{sec:Projection}

In this section we briefly discuss the construction of the projection operators $\RL$ and $\RR$ of degree $p$, which are applied in Sec. \ref{subsec:Coupling Degree}. Here we consider tensor-product based projection operators. As for the tensor-product SBP operators, we derive the projection operators in one dimension, denoted by $\RLone$ and $\RRone$.
We focus on a left element $L$ and a right element $R$ with nodal distributions $\bm{\eta}_L$ and $\bm{\eta}_R$. Let $\bm{\eta}_{\Gamma}$ denote the nodal distribution on the intermediate grid, then \eqref{eq:RLRR} and \eqref{eq:dofinal} give us
\begin{equation}\label{eq:AppRL}
\begin{split}
\RLone\etakL&=\etakGam,\\
\left(\RLone\right)\Tr\MGam\etakL&=\MetaLone\etakGam,
\end{split}
\end{equation}
and
\begin{equation}\label{eq:AppRR}
\begin{split}
\RRone\etakR&=\etakGam,\\
\left(\RRone\right)\Tr\MGam\etakR&=\MetaRone\etakGam,
\end{split}
\end{equation}
for $k=0,\dots,p$. We reiterate, that the choice of the nodes on the intermediate grid is nearly arbitrary, the only requirement is that a $2p$-accurate quadrature rule exists. In the case of $\bm{\eta}_{\Gamma}=\bm{\eta}_L$ or $\bm{\eta}_{\Gamma}=\bm{\eta}_R$, we set the projection operators to be the identity matrix.
If we consider the Gauss nodes on the intermediate grid as in Sec. \ref{subsec: Conservation and energy stability}, the projection matrices are fully determined in the sense that no free parameters remain.
In general, (\ref{eq:AppRL}) and (\ref{eq:AppRR}) are insufficient to fully specify $\RLone$ and $\RRone$. We want the approximation of the surface integral to be as close as possible to the approximated surface integral for the conforming case, in a sense that $\MetaLone\approx\left(\RLone\right)\Tr\MGam\RLone$ and $\MetaRone\approx\left(\RRone\right)\Tr\MGam\RRone$. Therefore, we use any remaining degrees of freedom to ensure that the modulus of each eigenvalue of the matrices $\MetaLone-\left(\RLone\right)\Tr\MGam\RLone$ and $\MetaRone-\left(\RRone\right)\Tr\MGam\RRone$ is as close to zero as possible.

Similarly, we motivate this optimization by analyzing the numerical flux in the left element for the upwind SAT, i.e. for $\sigma=1$.
For conforming nodal distributions, this numerical flux in one dimension for the linear advection equation, with unit wave speed, reduces to
\begin{equation}
\bm{f}^{*,(1\textrm{D})}_{L,conf}=\lamxi\RxiLN\Tr\Metaone\RxiLN\uL.
\end{equation}
For $\sigma=1$, the numerical flux of the non-conforming approximation \eqref{eq:ustarLR} reduces to
\begin{equation}
\bm{f}^{*,(1\textrm{D})}_{L}=\frac{\lamxi}{2}\RxiLN\Tr\left(\MetaLone+\left(\RLone\right)\Tr\MGam\RLone\right)\RxiLN\uL.
\end{equation}
Here, $\bm{f}^{*,(1\textrm{D})}_L=\bm{f}^{*,(1\textrm{D})}_{L,conf}$, if $\MetaLone=\left(\RLone\right)\Tr\MGam\RLone$.

Our approach to optimize the eigenvalues of $\MetaLone-\left(\RLone\right)\Tr\MGam\RLone$ and $\MetaRone-\left(\RRone\right)\Tr\MGam\RRone$ is based on the work of Kozdon and Wilcox~\cite{Kozdon2016}. Here, a Levenberg-Marquardt optimization algorithm has been used. We set the number of optimization iterations within the algorithm to $10$.
For constructing one-dimensional projection operators a MATLAB code is provided \textred{in the electronic supplementary material (ESM).}


\section{Comparison to the non-conforming method of Kozdon and Wilcox}\label{sec:KW}
To put the derivations using the multi-dimensional notation into context, it is useful to observe that for the tensor-product ansatz our approach is equivalent to the scheme of Kozdon and Wilcox \cite{Kozdon2016}. In their notation, a projection from the left element to the glue grid is denoted by the one-dimensional projection operator $\mat{P}_{L2G}$. The back transformation is denoted by $\mat{P}_{G2L}$. If we also consider a tensor-product approximation and take the glue grid to be equivalent to our intermediate grid, then
\begin{equation}
\mat{P}_{L2G}=\RLone.
\end{equation}
Due to the accuracy conditions made for $\RLone$, the operator $\mat{P}_{L2G}$ is also $p$-th degree.

For the approach of Kozdon and Wilcox, as for the approach of Mattsson and Carpenter \cite{Mattsson2010b}, the projection operators need to satisfy the following condition to ensure stability:
\begin{equation}\label{eq:StabKW}
\MetaLone\mat{P}_{G2L}=\mat{P}_{L2G}\Tr\MGam.
\end{equation}
Motivated by the stability condition one can construct a transformation back to left element
\begin{equation}\label{eq:KW}
\mat{P}_{G2L}=\left(\MetaLone\right)^{-1}\left(\RLone\right)\Tr\MGam^{(1\textrm{D})}.
\end{equation}
Multiplying \eqref{eq:KW} by a monomial, say $\etakGam$ with $k=0,\dots,p$, and using \eqref{eq:dofinal} we find
\begin{equation}
\mat{P}_{G2L}\etakGam=\left(\MetaLone\right)^{-1}\underbrace{\left(\RLone\right)\Tr\MGam \etakGam}_{=\MetaLone\etakL}=\etakL.
\end{equation}
So the $\mat{P}_{G2L}$ is also $p$-th accurate. Considering the symmetric SAT in the discretization (\ref{eq:discA} with $\sigma=0$) and focusing on tensor-product operators, we get
\begin{equation}\label{eq:KWSAT}
\bm{SAT}=\frac{\lamxiA}{2\fnc{J}_L}\HA^{-1}\left(\RxiLN\Tr\MetaLone\RxiLN\uA-\RxiLN\Tr\left(\RLone\right)\Tr\MGam\RRone\RxiRone\uR\right).
\end{equation}
Denoting $\uLN\coloneqq \RxiLN\uL$ and $\uRone\coloneqq \RxioneR\uR$ we can rewrite \eqref{eq:KWSAT} to become
\begin{equation}
\begin{split}
SAT&=\frac{\lamxiA}{2\fnc{J}_L}\left(\uLN-\underbrace{(\MetaLone)^{-1}(\RL^{(1\textrm{D})})\Tr\MGam}_{=\mat{P}_{G2L}}\underbrace{\RR^{(1\textrm{D})}}_{=\mat{P}_{R2G}}\uRone\right),\\
&=\frac{\lamxiA}{2\fnc{J}_L}\left(\uLN-\mat{P}_{G2L}\mat{P}_{R2G}\uRone\right),
\end{split}
\end{equation}
which is precisely the SAT for the Kozdon and Wilcox approach, assuming a symmetric SAT, despite the fact we have a different approach to derive the SAT.

\end{document}